# 2段階スケジューリング法を用いた看護師勤務表自動作成支援とその実運用

仲島 圭将　　古池 浩平　　井上 文彰
大阪大学

**和文概要**　　看護師勤務表の作成は看護の質・安全性，ならびに看護師の生活の質を左右する非常に重要な業務であり，日本のほとんどの病院では，病棟の管理者である看護師長がその責任を負って実施している．勤務表作成にかかる身体的・精神的負担は非常に大きいことに加え，近年における看護師不足の深刻化や働き方の多様化により作成の困難さが増していることから，看護師勤務表の自動作成に対する需要は年々高まっている．技術的には，現代の整数計画問題ソルバを用いれば現実的な時間で勤務表の自動作成が可能である．一方，未だ多くの病院において勤務表の作成は手作業で行われているという現実がある．これは，看護師長が無意識的に考慮していること（暗黙知）のすべてを制約条件として定式化できないため，最適化問題の解として出力された勤務表が実運用に耐えないものとなってしまうためである．本稿では，この問題を解決するための新たなアプローチとして，「2段階スケジューリング法」を提案する．この手法は，勤務表作成問題を夜勤スケジューリングと日勤スケジューリングに分割し，その途中段階で看護師長による手直しを行うことで，実際に現場で利用可能な看護師勤務表の作成を可能とするという発想に基づいている．また，看護師勤務表自動作成の普及の一助とすべく，実際の急性期病院および慢性期病院の現場において提案手法に基づくシステムを実運用した事例を紹介し，そこから得られた課題およびその対策について検討する．

**キーワード**: スケジューリング, 看護師勤務表, 2段階スケジューリング, 自動作成支援, 実運用

## 1. はじめに

看護師勤務表は看護の質・安全性と看護師の生活の質を左右する非常に重要な業務であるため，日本の多くの病院では，従来より，病棟の管理者である看護師長が勤務表作成の責任を負っている．看護師長が勤務表を手作業で作成する場合にかかる平均時間は10時間程度と言われており [9]，この作業が苦手な看護師長は自宅に持ち帰って作成することもある．また，勤務表作成の困難さを理解していないスタッフも少なくなく，看護師長が長時間をかけて多様な事情を勘案した勤務表を作成したとしても，スタッフが不平や不満を漏らすことも多い．このため，勤務表作成は，看護師長にとって身体的にも精神的にも負担が大きい業務と言える．さらに，少子高齢化に伴う看護師不足 [11] や働き方の多様化により，勤務表作成自体の難易度も上がっており，自動作成に対する需要は年々高まっている．

看護師勤務表の自動作成は，ナース・スケジューリング問題（Nurse Scheduling Problem; NSP）[2] と呼ばれ，看護の質や安全性を表すシフト制約（勤務表の縦の制約）と看護師一人一人の生活の質を表すナース制約（勤務表の横の制約）を満たすように，日勤，夜勤，休み等のシフトを割り当てる0-1整数計画問題として定式化される．日本におけるNSP研究は，1990年代後半の池上らの研究 [3] に端を発し，定式化やヒューリスティックアルゴリズムの開発に関する研究が中心であった．これは，当時の計算機の処理能力では，NSPの最適解を得ることが非常に難しかったためである．2010年ごろから，計算機の処理能力の向上に加え，GurobiやCPLEX等に代表される最適化ソルバの性能が向上したことで，実際の看護現場の条件を反映したNSPの最適解を得ることが可能となり，近年では，多くの看護師勤務表自動作成システムが商用化されている [7, 8, 14, 15]．

このように看護師勤務表を自動作成する技術的基盤は整っているにもかかわらず，未だ多くの病院では手作業での勤務表作成がなされているという現実がある．著者の仲島は，看護師長以上の役職者が集まる会議である看護管理学会学術集会において，2022年から継続的にNSPに関する発表を行い[18, 19, 20]，現場の看護管理者と意見交換を行ってきた．そこでは，「自動作成できるなんて知らなかった」，「自動作成システムを試験導入したが，思ったものができず，結局やめてしまった」，「自動作成システムを試験導入したが，うまく横展開できない」という現場の声が多く聞かれた．このように，技術的には十分に自動作成が可能であり，その需要も大きいにもかかわらず，現場レベルには普及していないというギャップが生じている背景には，自動作成の認知度の問題のほか，次の二つの課題があると考えられる．

**課題1.** 看護師長が手作業で作成する際に無意識的に考慮していること（暗黙知）をすべて制約条件として定式化することが難しい．例えば，看護師Aは，通常，夜勤においてリーダー的な役割を担うこと（1番手）が可能だけれども，その日の2番手が看護師Bの場合は難しい（ただし，3番手が看護師Cであれば可能）というような「微妙な人間関係」に関する制約である．よって，必然的に自動作成された勤務表を修正することになるが，非常に多くの条件を考慮して自動作成された「最適な」勤務表を部分的に手修正することは一般に困難である．

**課題2.** 実際の看護現場では，採用や離職を含む人事異動や人間関係の変化が頻繁に発生するため，基本的には，スケジューリングを実行する都度(すなわち毎月)設定を見直す必要があるが，数理最適化に関する知識がほとんどない看護師長には「どの設定を変更するとどの制約条件に反映されるのか」を把握することが難しく，その見直し作業自体が困難である．ある程度慣れて使いこなせるようになったとしても，その知識は看護師長個人に依存してしまうため，看護師長が異動や退職した場合にはその知識が失われる可能性が高く，継続的な運用が難しい．

本稿の目的は，上記の課題を解決した事例を紹介することである．まず，上記の課題1を解決するための新たなアプローチとして，**2段階スケジューリング法**を提案する．多くの商用システムでは，設定された制約条件などをもとに，日勤，夜勤，ならびに休みを含むすべての勤務記号を一括して配置する最適化問題を考え，その求解を行うというアプローチが取られている．これに対し，本稿で提案する2段階スケジューリング法では，夜勤のみを配置する**夜勤スケジューリング**と日勤（早出，遅出を含む）および休みを配置する**日勤スケジューリング**の2段階に勤務表作成問題を分割し，初めに夜勤スケジューリングのみを実行した段階で**ユーザ(看護師長)による手直し**を実施する．さらに，それに続いて日勤スケジューリングを実行することで最終的な勤務表を得る．これまで，一つの最適化問題として定式化し，それを夜勤，日勤のサブ問題に分割して，夜勤，日勤の順番に解くヒューリスティックな手法[1]は存在したが，本稿で提案する2段階スケジューリング法は，夜勤スケジューリングと日勤スケジューリングをそれぞれ別問題として定式化するという点で新しいアプローチである．このように完全に問題を分割することで，夜勤スケジューリングを解いた時点では，まだ空欄（勤務が配置されていない箇所）が多いため，修正自体が容易になる．さらに，一般に，日勤に比べて夜勤の配置人数は少ないため，安全性の点でその重要度が高い．その夜勤を看護師長の考えに沿って容易に修正できるようにすることで，看護師長の納得感も高まり，実運用に耐える勤務表を最終的に生成することが可能となる．

加えて，本稿では，上記の課題2への対処に関する知見を提供することを目的として，著

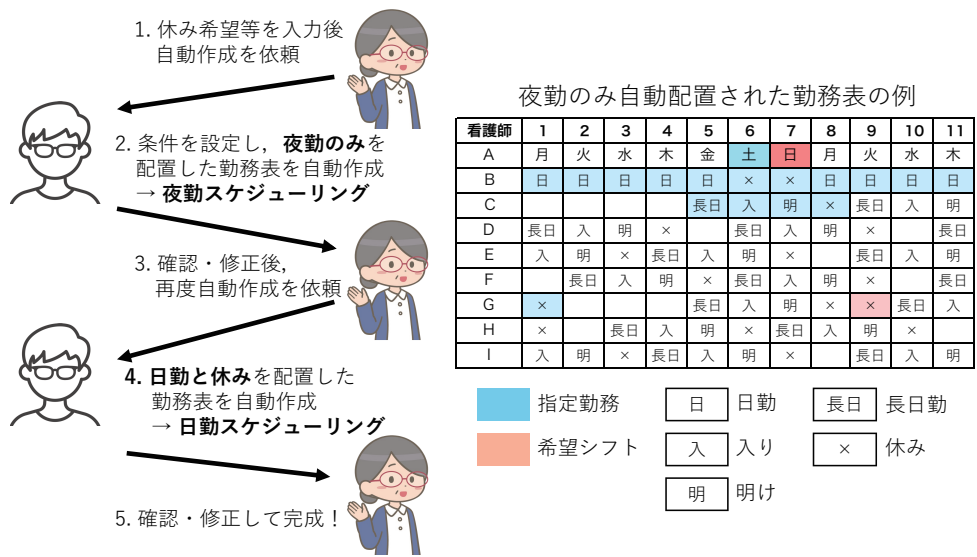

図 1: 2 段階スケジューリング法を用いた自動作成支援の流れ

者の仲島が開発ならびに運用している勤務表自動作成支援システム [13] を紹介する．特に，ケーススタディとして，このシステムを急性期病院である大阪大学医学部附属病院，ならびに慢性期病院である阪神リハビリテーション病院において実運用した事例について報告し，そこから得られた課題およびその対策について検討する．

以下の本稿は次のように構成される．まず，2 節において，2 段階スケジューリング法の概要および 2 段階スケジューリング法を用いた自動作成支援について説明する．次に，3 節において，大阪大学医学部附属病院の事例に基づき，夜勤スケジューリングおよび日勤スケジューリングの定式化を示す．4 節では，いくつかの具体事例を用いて，夜勤スケジューリングおよび日勤スケジューリングの実行時間を評価し，また，看護師長の手直しによる目的関数値の変化について述べる．また，5 節では急性期病院および慢性期病院において，2 段階スケジューリング法を用いた看護勤務表自動作成支援を実運用した事例とそこから得られた課題と対策について述べる．最後に，6 節において本稿を締めくくる．

なお本稿は，日本オペレーションズリサーチ学会 2022 年春季研究発表会での口頭発表 [17] に基づいている．

## 2. 2 段階スケジューリング法

### 2.1. 概要

前述の通り，本手法では，勤務表自動作成を夜勤スケジューリングと日勤スケジューリングの 2 段階に分割し，その途中で看護師長による手直しを実施する（図 1，勤務記号については 3 節を参照）．これは，「手作業で勤務表を作成する際，多くの看護師長は最初に夜勤を配置してから日勤と休みを配置する」という知見を取り入れたものである．このように先に夜勤を配置するのは，7 人の患者に対して 1 人の看護師が配置される 7 対 1 看護や 10 人の患者に対して 1 人の看護師が配置される 10 対 1 看護と呼ばれる標準的な看護体制の下では，夜勤人数は日勤人数の半分以下となり，かつ，夜勤帯は基本的に看護師長が不在となるため，安全性および看護師長の心理面から夜勤配置に重点が置かれることに起因している．

本手法では，最初に夜勤のみを配置した勤務表を自動作成し，その段階での手直しを許すことで，空欄（勤務が配置されていない箇所）が多く自由度の高い状態で，看護師長の意図

に沿った勤務表へと修正することが可能となる．特に，重要度の高い夜勤での人員配置に対して確認および修正ができるため，看護師長の納得感を高めることができる [18]．

なお，数理最適化の観点では，問題を 2 段階に分割し，かつ手直しを許していることで，最終的に得られる勤務表は明らかに原問題の最適解とは異なっている．しかし，このような手直しは実運用において非常に重要であり，この意味で，看護師勤務表の作成問題では定式化された問題の最適解は必ずしも運用上の良好な解とはならないことを示している．

### 2.2. 2 段階スケジューリング法を用いた自動作成支援

多くの既存サービスでは，勤務表自動作成システムのみを納品し，基本的な運用は各病院に委ねるという形態をとっている [7, 8, 14, 15]．FAQ やチャット等を用いた問い合わせ対応はあるが，基本的には看護師長自身で問題を解決する必要があるため，前節で述べた課題 2 に対して継続的に対応することは難しい．そこで，本手法に基づく実運用においては，以下の手順によって自動作成支援を行っている（図 1）．

1. 各病棟の看護師長は休み希望や会議・研修等の指定勤務を入力し，対象月の勤務希望表を作成する．また，対象月の設定条件をまとめた「勤務条件設定ファイル[1]」を作成し，前月分，対象月（および，必要であれば来月分）の勤務表と合わせてシステム担当者[2]に共有する．この共有過程については，看護部長等が全病棟分をまとめて管理し，一括で共有することも可能である．一括する方が人事異動等の管理面でメリットがあるが，この点については 5.2 節で述べる．

2. 勤務条件設定ファイルに基づき，担当者が条件を正しく反映させ，夜勤のみを配置した勤務表を自動作成する（夜勤スケジューリング）．その際，勤務希望において「明け–日勤」等の配置してはいけないシフトに相当するものが入力されている場合は，その内容を看護師長に確認し，修正の上，再度自動作成する．そして，自動作成された勤務表に加え，その評価（例えば，夜勤者数が明らかに不足している）を，フィードバックコメントとして報告する．

3. 手順 2 において作成された勤務表およびコメントをもとに，各看護師長が確認・修正を行い，その修正したファイルを担当者に共有する．必要に応じて，夜勤スケジューリングを再度行う場合もある．

4. 手順 3 において修正された勤務表をもとに，担当者が日勤（早出と遅出を含む）と休みを配置した勤務表を自動作成する（日勤スケジューリング）．手順 2 と同様に，自動作成した勤務表とフィードバックコメントを共有する．

5. 手順 4 において作成された勤務表およびコメントをもとに，各看護師長が確認・修正を行い，勤怠管理システムに反映させる．

このように，看護師長が（現場の文脈で）作成した条件設定ファイルに基づいて，担当者がその条件を（数理最適化の文脈に置き換えて）正しく自動作成システムに反映させることで課題 2 に対応する．看護師長が作成する条件設定ファイルには，システム上設定不可能な条件が含まれることもあるが，その場合は，担当者が可能な範囲の条件設定を行い，完全に実現することは不可能な旨も含めて適切にフィードバックする．

---

[1] 制約条件からフォーマットを決定．数値や記号のみの入力および自由記載による説明を含む．
[2] 後述の事例では著者の仲島が担当．

## 2.3. 自動作成支援のための実装上の工夫

自動作成支援を効率的に行うためには，看護師長が作成する勤務希望表や条件設定ファイルにおいて，求解の障害となる箇所を素早く発見することが重要である．本手法では，実装上，「勤務希望表におけるマスタに存在しない勤務記号の検出」と「条件設定ファイルにおける未設定箇所の検出」に加え，「ハード制約チェック」の三つの工夫を行っている．ここでは，特に，看護師長による勤務表の修正のしやすさにも関わる「ハード制約チェック」について詳しく説明する．

一般に，制約条件は，必ず満たす必要があるハード制約と，できるだけ満たすことが望ましいソフト制約に分類される．ハード制約を多く設定することで，看護師長が意図した通りの勤務表が作成される可能性が高くなるが，例外的な事情によりハード制約を満たすことが困難な場合，最適解が得られない．例えば，「夜勤明けの後は休みにする」というハード制約を設定している場合，たとえ本人が許容した上で夜勤明けの後に研修を配置していたとしても，このハード制約を満たさないため，最適解が得られない．このような事態を防ぐために，既存システムの中には，ハード制約をほとんど設定しないものもあるが，その場合，看護師長が意図した通りの勤務表が作成される可能性が低くなり，結果的に多大な修正が必要となる．そこで，本システムでは，ハード制約を多く設定し，看護師長が意図した通りの勤務表が作成される可能性を高めるとともに，夜勤スケジューリングおよび日勤スケジューリングを解く前に，ハード制約を緩和した問題を解くことで違反しているハード制約がないかを確認する．もし満たされていないハード制約がある場合には，そのハード制約から原因となっている箇所を特定することで，すぐに看護師長に確認を取ることが可能となる．

## 3. 2段階スケジューリングの定式化

本節では，大阪大学医学部附属病院の事例に基づき，2段階スケジューリングの0-1整数計画問題としての定式化を示す．これまで大阪大学医学部附属病院では，「夜勤入り–夜勤明け–休み」が基本配置である2交代制16時間夜勤を採用してきたが，2022年度から2交代制12時間夜勤および12時間日勤（通常の8時間勤務の日勤に対し，長日勤とも呼ばれる）に移行しており，「12時間日勤–夜勤入り–夜勤明け–休み」が基本配置となっている [10]．

本稿において利用する全シフト集合 $\mathcal{S}$ を以下として与える．

$$\mathcal{S} = \{\, 日, 12h, 早, 遅, 入, 明, 他, 休, 特休, 未 \,\}$$

ただし，「日」，「12h」は通常の（8時間）日勤および12時間日勤を表し，「早」，「遅」はそれぞれ早出，遅出を表す．「入」，「明」は夜勤入りおよび夜勤明けを表し，「他」は会議や研修等の病棟以外での日勤勤務を表す．「休」は土日および祝日（振替も含む）に相当する休み（週休とする）を表し，「特休」は有給休暇や夏季休暇等の週休以外の休み（特別休暇）を表す．「未」は未確定を表し，夜勤スケジューリングにおいてのみ利用される．勤務表において実際に利用される勤務記号は100種類以上存在するが，本システムではこれらを上記の9種類に集約している．

ここで，定式化において用いる集合および変数に関する記号を表1にまとめる．看護師には，日勤と12時間夜勤の両方に従事する看護師（夜勤可能看護師とする），16時間夜勤のみに従事する夜勤専従看護師，日勤のみに従事する日勤専従看護師（時短看護師も含む）が存在する．また，基本的には患者を受け持たず，病棟全体の動きが円滑になるように調整するリーダー看護師（日勤リーダー）が存在する．次に，グループは看護実践能力別の看護

師集合を表す．本稿では 4 つのグループを考え，番号が小さいほど実践能力が高いものとする．チームはグループを跨ぐ縦割り集合であり，本稿では，$A$ と $B$ の 2 チームを考える．疾患別チームやフロア別チーム等の用途が考えられるが，病棟によってはチーム分けを行わない場合もある．ここで，各看護師はグループまたはチームに必ず所属するわけではないことに注意する．例えば，大阪大学医学部附属病院の場合，看護師長はグループにもチームにも所属しないことがほとんどである．表中の「曜日依存」とある変数については，実際には曜日ごとに値を設定しているが，本稿では簡単のために定数としていることを表す．

表 1: 集合および変数一覧

| 記号 | 説明 |
| --- | --- |
| $\mathcal{N}$ | 全看護師集合 |
| $\mathcal{N}_\text{n} \subset \mathcal{N}$ | （夜勤専従除く）夜勤可能看護師集合 |
| $\mathcal{N}_\text{no} \subset \mathcal{N}$ | 夜勤専従看護師集合 |
| $\mathcal{N}_\text{day} \subset \mathcal{N}$ | 時短を含む日勤専従看護師集合 |
| $\mathcal{N}_\text{r} \subset \mathcal{N}$ | 新人看護師（rookie）集合 |
| $\mathcal{N}_\text{men} \subset \mathcal{N}$ | 男性看護師集合 |
| $\mathcal{N}_\text{dayL} \subset \mathcal{N}$ | 日勤リーダーが可能な看護師集合 |
| $\mathcal{S}$ | 全シフト集合 |
| $\mathcal{S}_\text{off}, \mathcal{S}'_\text{off}$ | 休みシフト集合：$\mathcal{S}_\text{off} = \{$ 休,特休 $\}$，$\mathcal{S}'_\text{off} = \mathcal{S}_\text{off} \cup \{$ 未 $\}$ |
| $\mathcal{S}_\text{work}, \mathcal{S}'_\text{work}$ | 勤務シフト集合：$\mathcal{S}_\text{work} = \{$ 日,12h,早,遅,入,明,他 $\}$，$\mathcal{S}'_\text{work} = \mathcal{S}_\text{work} \cup \{$ 未 $\}$ |
| $\mathcal{S}_\text{day}$ | 日勤帯に関するシフト集合：$\{$ 日,12h,早,遅,他 $\}$ |
| $\mathcal{S}_\text{hf}, \mathcal{S}_\text{sf}$ | ハード制約およびソフト制約における禁止シフト集合 |
| $\boldsymbol{s}^{(\text{h})}, \boldsymbol{s}^{(\text{s})}$ | $\mathcal{S}_\text{hf}, \mathcal{S}_\text{sf}$ の要素であり，$\{$ 明,入 $\}$ のような禁止するシフト系列の順序付き集合を表す |
| $\mathcal{G}$ | グループ集合：$\{1,2,3,4\}$ |
| $\mathcal{N}_g \subset \mathcal{N}$ | グループ $g\ (\in \mathcal{G})$ に属する看護師集合 |
| $\mathcal{N}_{12}, \mathcal{N}_{34} \subset \mathcal{N}$ | $\mathcal{N}_{12} = \mathcal{N}_1 \cup \mathcal{N}_2, \mathcal{N}_{34} = \mathcal{N}_3 \cup \mathcal{N}_4$ |
| $\mathcal{N}_\text{grp} \subset \mathcal{N}$ | グループに属する看護師集合：$\bigcup_{g \in \mathcal{G}} \mathcal{N}_g$ |
| $\mathcal{T}$ | チーム集合：$\{A, B\}$ |
| $\mathcal{N}_t \subset \mathcal{N}$ | チーム $t\ (\in \mathcal{T})$ に属する看護師集合 |
| $\mathcal{N}_{t,g} \subset \mathcal{N}$ | $\mathcal{N}_t \cap \mathcal{N}_g$ |
| $\mathcal{P}$ | ペア集合：$\{p_i = (n_1, n_2, s_1, s_2, N_{i,\min});\ i = 1, 2, \ldots, |\mathcal{P}|,\ n_1, n_2 \in \mathcal{N},\ s_1, s_2 \in \mathcal{S}_\text{work}\}$，$N_{i,\min}$ はペア $p_i$ における同日勤務回数の下限 |
| $\mathcal{P}_\text{f}$ | 禁止ペア集合：$\{p_i^{(\text{f})} = (n_1, n_2, s_1, s_2);\ i = 1, 2, \ldots, |\mathcal{P}_\text{f}|,\ n_1, n_2 \in \mathcal{N},\ s_1, s_2 \in \mathcal{S}_\text{work}\}$ |
| $\mathcal{F}$ | 禁止勤務に関する集合：$\{f_i = (n, s, d_\text{w});\ i = 1, 2, \ldots, |\mathcal{F}|,\ n \in \mathcal{N},\ s \in \mathcal{S}_\text{work}\}$，$d_\text{w}$ は曜日を表す |
| $\mathcal{D}_\text{p}, \mathcal{D}_\text{n}$ | 前月末/次月頭の 5 日間の日付集合 |
| $\mathcal{D}_\text{t}$ | 勤務表を作成する対象月の日付集合 |
| $\mathcal{D}_\text{t,off}$ | 勤務表を作成する対象月の土日祝にあたる日付集合 |
| $\mathcal{D}$ | $\mathcal{D}_\text{p} \cup \mathcal{D}_\text{t} \cup \mathcal{D}_\text{n}$ |

| 記号 | 説明 |
|---|---|
| $\mathcal{D}_{5d} \subset \mathcal{D}$ | 夜勤を含む「5 連勤–1 休み–4 連勤」または「4 連勤–1 休み–5 連勤」を回避する制約のための日付集合 |
| $\mathcal{D}_{5d}^{(d)} \subset \mathcal{D}$ | （夜勤を含む）5 連勤を回避する制約のための日付集合 |
| $\mathcal{D}_{5d1off} \subset \mathcal{D}$ | 5 連続勤務の後の 1 休みを回避する制約のための日付集合 |
| $\mathcal{D}_{9d2off} \subset \mathcal{D}$ | 9 日間における週休回数の下限に関する制約のための日付集合 |
| $\mathcal{D}_{6d} \subset \mathcal{D}$ | 6 連続勤務を回避する制約のための日付集合 |
| $\mathcal{D}_{3d} \subset \mathcal{D}$ | 12 時間日勤の前の 3 連続勤務を回避する制約のための日付集合 |
| $\mathcal{D}_{3d}^{(d)} \subset \mathcal{D}$ | 3 連続日勤を回避する制約のための日付集合 |
| $\mathcal{D}_{3n} \subset \mathcal{D}$ | 3 連続夜勤を回避する制約のための日付集合 |
| $\mathcal{D}_{4n}, \mathcal{D}_{4no} \subset \mathcal{D}$ | 4 連続の 12 時間/16 時間夜勤を回避する制約のための日付集合 |
| $\mathcal{D}_{2w4n} \subset \mathcal{D}$ | 2 週間で 4 回以上の夜勤を回避する制約のための日付集合 |
| $\mathcal{D}_{iub}, \mathcal{D}_{ilb} \subset \mathcal{D}$ | 夜勤インターバルの上下限に関する制約のための日付集合 |
| $\mathcal{D}_s \subset \mathcal{D}$ | シフト $s$ ($\in \mathcal{S}_{work}$) に対する同一勤務の連続回数制限のための日付集合 |
| $N_{s,max}$ | シフト $s$ ($\in \mathcal{S}_{work}$) の最大連続回数 |
| $\mathcal{D}_{3off}$ | 3 連休候補を 2 箇所以上確保する制約のための日付集合 |
| $\mathcal{D}_{2off}^{(d)}$ | 2 連休回数の下限制約のための日付集合 |
| $\mathcal{D}_{3off}^{(d)}$ | 3 連続以上の休み回数の下限制約のための日付集合 |
| $\mathcal{D}_f \subset \mathcal{D}$ | 禁止シフト制約のための日付集合 |
| $\mathcal{D}_{fssm} \subset \mathcal{D}_t$ | 対象月の金土日月の日付集合 |
| $\mathcal{D}_{sat} \subset \mathcal{D}_t$ | 対象月の土曜日の日付集合（最終日の土曜日は含まない） |
| $\mathcal{D}_{ope} \subset \mathcal{D}_t$ | 対象月の手術日等を含む指定した曜日の日付集合 |
| $d_{p,last} \in \mathcal{D}_p$ | 前月の最終日 |
| $d_1, d_{last} \in \mathcal{D}_t$ | 対象月の 1 日/最終日 |
| $d_{n1}, d_{n2} \in \mathcal{D}_n$ | 次月の 1 日/2 日 |
| $f_n^{(n)} \in \{0,1\}$ | 看護師 $n$ の夜勤回数制限フラグ：$f_n^{(n)} = 1$ の場合，指定回数で固定する |
| $f_n^{(p)}$ | 看護師 $n$ の休みの好み：$f_n^{(p)} \in \{$ 土日連休,3 連休,連休,単発休み $\}$ |
| $\mathcal{N}_{so} \subset \mathcal{N}$ | 単発休みが好みの看護師集合：$\{n;\ n \in \mathcal{N},\ f_n^{(p)} = $ 単発休み $\}$ |
| $x_{n,d,s} \in \{0,1\}$ | 看護師 $n$，日付 $d$，シフト $s$ に対する決定変数 |
| $\boldsymbol{X}$ | 決定変数を並べたベクトル |
| $x'_{n,d,s} \in \{0,1\}$ | 看護師 $n$，日付 $d$，シフト $s$ に対する希望または指定勤務を表す変数 |
| $N_{off}$ | 対象月に配置可能な休み（週休，祝日，振替）回数 |
| $N_n^{(UB)}, N_n^{(LB)}$ | 看護師 $n$ の夜勤回数の上下限 |
| $\bar{N}_{fssm}$ | 対象月の土日回数および夜勤可能看護師数に基づいて算出された金土日月の平均夜勤回数 |
| $\bar{N}_{ope}$ | 対象月における手術等のイベント日および夜勤可能看護師数に基づいて算出されたイベント日の平均夜勤回数 |
| $w \in \{d, n\}$ | 勤務帯（「日, 12h, 早, 遅」を含む日勤帯：d, 夜勤帯：n） |
| $N_w^{(UB)}, N_w^{(LB)}$ | $w$ に配置する看護師数の上下限（曜日依存） |
| $N_{w,1}^{(UB)}, N_{w,1}^{(LB)}$ | $w$ に配置するグループ 1 の勤務人数の上下限（曜日依存） |
| $N_{w,12}^{(UB)}, N_{w,12}^{(LB)}$ | $w$ に配置するグループ 1, 2 の勤務人数の上下限（曜日依存） |

| 記号 | 説明 |
| --- | --- |
| $N_{w,4}^{(\mathrm{UB})}, N_{w,4}^{(\mathrm{LB})}$ | $w$ に配置するグループ4の勤務人数の上下限（曜日依存） |
| $N_{w,t,1}^{(\mathrm{UB})}, N_{w,t,1}^{(\mathrm{LB})}$ | $w$ に配置するチーム$t$かつグループ1の勤務人数の下限（曜日依存） |
| $N_{w,t,12}^{(\mathrm{UB})}, N_{w,t,12}^{(\mathrm{LB})}$ | $w$ に配置するチーム$t$かつグループ1，2の勤務人数の上下限（曜日依存） |
| $N_{w,t,34}^{(\mathrm{UB})}, N_{w,t,34}^{(\mathrm{LB})}$ | $w$ に配置するチーム$t$かつグループ3，4の勤務人数の上下限（曜日依存） |
| $N_{w,t}^{(\mathrm{UB})}, N_{w,t}^{(\mathrm{LB})}$ | $w$ に配置するチーム$t$の勤務人数の上下限（曜日依存） |
| $N_{w,\mathrm{r}}^{(\mathrm{UB})}, N_{w,\mathrm{r}}^{(\mathrm{LB})}$ | $w$ に配置する新人看護師の勤務人数の上下限（曜日依存） |
| $N_{w,\mathrm{men}}^{(\mathrm{UB})}, N_{w,\mathrm{men}}^{(\mathrm{LB})}$ | $w$ に配置する男性看護師の勤務人数の上下限（曜日依存） |

表 2: 夜勤スケジューリングにおけるハード制約一覧

| 識別子 | 集合 | 内容 |
| --- | --- | --- |
| $\mathrm{H_n}$-N-1 | $\mathcal{H}_{\mathrm{n},1}$ | 勤務希望と指定勤務の固定および夜勤シフトのみに限定 |
| $\mathrm{H_n}$-N-2 | $\mathcal{H}_{\mathrm{n},2}$ | 前月からの月跨ぎの禁止シフト |
| $\mathrm{H_n}$-N-3 | $\mathcal{H}_{\mathrm{n},3}$ | 前月からの月跨ぎの夜勤シフト固定 |
| $\mathrm{H_n}$-N-4 | $\mathcal{H}_{\mathrm{n},4}$ | 次月への月跨ぎの固定・禁止シフト |
| $\mathrm{H_n}$-N-5 | $\mathcal{H}_{\mathrm{n},5}$ | 6連勤の禁止 |
| $\mathrm{H_n}$-N-6 | $\mathcal{H}_{\mathrm{n},6}$ | 禁止シフト |
| $\mathrm{H_n}$-N-7 | $\mathcal{H}_{\mathrm{n},7}$ | 日勤専従看護師の夜勤制限 |
| $\mathrm{H_n}$-N-8 | $\mathcal{H}_{\mathrm{n},8}$ | 夜勤可能および夜勤専従看護師の夜勤回数制限 |
| $\mathrm{H_n}$-N-9 | $\mathcal{H}_{\mathrm{n},9}$ | 夜勤インターバルの下限 |
| $\mathrm{H_n}$-N-10 | $\mathcal{H}_{\mathrm{n},10}$ | 同一勤務の連続回数制限 |
| $\mathrm{H_n}$-N-11 | $\mathcal{H}_{\mathrm{n},11}$ | 3連休候補を2箇所以上確保 |
| $\mathrm{H_n}$-S-12 | $\mathcal{H}_{\mathrm{n},12}$ | 日勤においてリーダー業務可能な看護師を一人以上確保 |
| $\mathrm{H_n}$-S-13 | $\mathcal{H}_{\mathrm{n},13}$ | 男性看護師の人数制限 |
| $\mathrm{H_n}$-S-14 | $\mathcal{H}_{\mathrm{n},14}$ | 必ず勤務記号が割り当てられる |

### 3.1. 夜勤スケジューリング

夜勤スケジューリングにおいては，「入」，「明」，「休」，「未」の4つのシフトのみを配置し，その他のシフトは配置しないようにする．また，「12h」については，後工程として「入」の前日に配置する．看護師 $n \in \mathcal{N}$ の日付 $d \in \mathcal{D}$ にシフト $s \in \mathcal{S}$ が割り当てられた場合に 1，そうでない場合に 0 をとる決定変数 $x_{n,d,s} \in \{0,1\}$ を導入し，すべての決定変数を並べたベクトルを $\boldsymbol{X}$ と表す．

夜勤スケジューリングにおけるハード制約を表す集合族 $(\mathcal{H}_{\mathrm{n},j})_{j=1,2,\ldots,N_{\mathrm{n}}^{(\mathrm{h})}}$ ($\mathcal{H}_{\mathrm{n},j} \subseteq \{0,1\}^{|\mathcal{N}| \times |\mathcal{D}| \times |\mathcal{S}|}$) ならびにソフト制約違反によるペナルティを表す関数 $f_{\mathrm{n},i}(\boldsymbol{X})$ ($i = 1, 2, \ldots, N_{\mathrm{n}}^{(\mathrm{s})}$) を用いて，池上らの定式化 [2] と同様に，夜勤スケジューリング問題は以下で表現される．

$$\text{minimize} \quad \sum_{i=1}^{N_{\mathrm{n}}^{(\mathrm{s})}} \alpha_{\mathrm{n},i} f_{\mathrm{n},i}(\boldsymbol{X}) \quad \text{s.t.} \quad \boldsymbol{X} \in \bigcap_{j=1}^{N_{\mathrm{n}}^{(\mathrm{h})}} \mathcal{H}_{\mathrm{n},j}$$

ここで，$N_n^{(h)}$ および $N_n^{(s)}$ は夜勤スケジューリングにおけるハード制約およびソフト制約の総数を表し，$\alpha_{n,i}$ は $i$ 番目のソフト制約に対する重みを表す．ただし，重み $\alpha_{n,i}$ については，商用化サービス [12, 13] との関係上，本稿では公開しない．

夜勤スケジューリングにおけるハード制約一覧を表2に，ソフト制約一覧を表3に示す[3,4]．ただし，識別子における2番目のNはナース制約，Sはシフト制約を表す．各制約の詳細は付録Aならびに付録Bに記載する．

表 3: 夜勤スケジューリングにおけるソフト制約一覧

| 識別子 | 関数 | 内容 |
|---|---|---|
| $S_n$-S-1, 2 | $f_{n,1}, f_{n,2}$ | 日毎の勤務人数の下限・上限 |
| $S_n$-S-3, 4 | $f_{n,3}, f_{n,4}$ | グループ1の勤務人数の下限・上限 |
| $S_n$-S-5, 6 | $f_{n,5}, f_{n,6}$ | グループ1，2の勤務人数の下限・上限 |
| $S_n$-S-7, 8 | $f_{n,7}, f_{n,8}$ | グループ4の勤務人数の下限・上限 |
| $S_n$-S-9, 10 | $f_{n,9}, f_{n,9}$ | チーム $t$ かつグループ1の勤務人数の下限・上限 |
| $S_n$-S-11, 12 | $f_{n,11}, f_{n,12}$ | チーム $t$ かつグループ1，2の勤務人数の下限・上限 |
| $S_n$-S-13, 14 | $f_{n,13}, f_{n,14}$ | チーム $t$ かつグループ3，4の勤務人数の下限・上限 |
| $S_n$-S-15, 16 | $f_{n,15}, f_{n,16}$ | チーム $t$ の勤務人数の下限・上限 |
| $S_n$-S-17, 18 | $f_{n,17}, f_{n,18}$ | 新人看護師の勤務人数の下限・上限 |
| $S_n$-S-19 | $f_{n,19}$ | ペア勤務回数の下限 |
| $S_n$-S-20 | $f_{n,20}$ | 禁止ペア |
| $S_n$-N-21 | $f_{n,21}$ | 禁止勤務 |
| $S_n$-N-22 | $f_{n,22}$ | ソフト禁止シフト |
| $S_n$-N-23 | $f_{n,23}$ | 休み（週休，祝日，振替）回数を4回以上残す |
| $S_n$-N-24 | $f_{n,24}$ | 夜勤インターバルの上限 |
| $S_n$-N-25, 26 | $f_{n,25}, f_{n,26}$ | 夜勤回数の下限・上限 |
| $S_n$-N-27 | $f_{n,27}$ | 金土日月の夜勤回数の公平化 |
| $S_n$-N-28 | $f_{n,28}$ | オペ日等の夜勤回数の公平化 |
| $S_n$-N-29 | $f_{n,29}$ | 夜勤を含む5連続勤務の回避 |
| $S_n$-N-30 | $f_{n,30}$ | 夜勤を含む「5連勤–1休み–4連勤」および「4連勤–1休み–5連勤」の回避 |
| $S_n$-N-31 | $f_{n,31}$ | 3連続夜勤の回避 |
| $S_n$-N-32 | $f_{n,32}$ | 4連続夜勤の回避 |
| $S_n$-N-33 | $f_{n,33}$ | 2週間で4回以上の夜勤の回避 |
| $S_n$-N-34 | $f_{n,34}$ | 夜勤回数に関する違反度の最大最小化 |
| $S_n$-N-35 | $f_{n,35}$ | ナース制約に関する違反度の最大最小化 |

### 3.2. 日勤スケジューリング

決定変数 $x_{n,d,s}$ および希望または指定勤務を表す変数 $x'_{n,d,s}$ について，簡単のため，夜勤スケジューリングと同様の記法を用いる．目的関数については，夜勤スケジューリングと同

---

[3] $H_n$-S-13 は，対象病棟ごとに利用要否を選択可能としている．
[4] $S_n$-S-9 から $S_n$-S-16 のチームに関する制約および $S_n$-S-17, $S_n$-S-18 の新人に関する制約は，対象病棟ごとに利用要否を選択可能としている．

様の形式であるため省略する．また，日勤スケジューリングでは，夜勤スケジューリングの解を看護師長が修正したものを勤務希望とみなして $x'_{n,d,s}$ を固定する．

日勤スケジューリングにおけるハード制約一覧を表 4 に，ソフト制約一覧を表 5 に示す[5,6]．各制約の詳細は付録 C ならびに付録 D に記載する．

表 4: 日勤スケジューリングにおけるハード制約一覧

| 識別子 | 集合 | 内容 |
| --- | --- | --- |
| $H_d$-N-1 | $\mathcal{H}_{d,1}$ | 勤務希望の固定および日勤/休みシフトのみに限定 |
| $H_d$-N-2 | $\mathcal{H}_{d,2}$ | 夜勤専従看護師の日勤制限 |
| $H_d$-N-3 | $\mathcal{H}_{d,3}$ | 同一勤務の連続回数制限 |
| $H_d$-N-4 | $\mathcal{H}_{d,4}$ | 6 連勤の回避 |
| $H_d$-N-5 | $\mathcal{H}_{d,5}$ | 次月跨ぎの夜勤を含む 6 連勤の回避 |
| $H_d$-N-6 | $\mathcal{H}_{d,6}$ | 12 時間日勤の前の 3 連続日勤の回避 |
| $H_d$-N-7 | $\mathcal{H}_{d,7}$ | 禁止シフト |
| $H_d$-S-8 | $\mathcal{H}_{d,8}$ | 男性看護師の人数制限 |
| $H_d$-S-9 | $\mathcal{H}_{d,9}$ | 必ず勤務記号が割り当てられる |

表 5: 日勤スケジューリングにおけるソフト制約一覧

| 識別子 | 関数 | 内容 |
| --- | --- | --- |
| $S_d$-S-1, 2 | $f_{d,1}, f_{d,2}$ | 日毎の勤務人数の下限・上限 |
| $S_d$-S-3, 4 | $f_{d,3}, f_{d,4}$ | グループ 1 に属する看護師の勤務人数の下限・上限 |
| $S_d$-S-5, | $f_{d,5}, f_{d,6}$ | グループ 1，2 に属する看護師の勤務人数の下限・上限 |
| $S_d$-S-7, 8 | $f_{d,7}, f_{d,8}$ | グループ 4 に属する看護師の勤務人数の下限・上限 |
| $S_d$-S-9, 10 | $f_{d,9}, f_{d,10}$ | チーム $t$ かつグループ 1 の勤務人数の下限・上限 |
| $S_d$-S-11, 12 | $f_{d,11}, f_{d,12}$ | チーム $t$ かつグループ 1，2 の勤務人数の下限・上限 |
| $S_d$-S-13, 14 | $f_{d,13}, f_{d,14}$ | チーム $t$ かつグループ 3，4 の勤務人数の下限・上限 |
| $S_d$-S-15, 16 | $f_{d,15}, f_{d,16}$ | チーム $t$ の勤務人数の下限・上限 |
| $S_d$-S-17, 18 | $f_{d,17}, f_{d,18}$ | 新人看護師の勤務人数の下限・上限 |
| $S_d$-S-19 | $f_{d,19}$ | ペア勤務回数の下限 |
| $S_d$-S-20 | $f_{d,20}$ | 禁止ペア |
| $S_d$-S-21 | $f_{d,21}$ | リーダー業務可能な看護師の勤務人数の下限 |
| $S_d$-N-22 | $f_{d,22}$ | 禁止勤務 |
| $S_d$-N-23 | $f_{d,23}$ | ソフト禁止シフト |
| $S_d$-N-24, 25 | $f_{d,24}, f_{d,25}$ | 週休回数の下限・上限 |
| $S_d$-N-26 | $f_{d,26}$ | 土日祝における休み回数の下限 |
| $S_d$-N-27 | $f_{d,27}$ | 3 連続日勤の回避 |
| $S_d$-N-28 | $f_{d,28}$ | 5 連勤の後の単発休みの回避 |
| $S_d$-N-29 | $f_{d,29}$ | 9 日間における週休回数の下限 |

---

[5] $H_d$-S-8 は，対象病棟ごとに利用要否を選択可能としている．
[6] $S_d$-S-9 から $S_d$-S-16 のチームに関する制約および $S_d$-S-17, $S_d$-S-18 の新人に関する制約は，対象病棟ごとに利用要否を選択可能としている．

| 識別子 | 関数 | 内容 |
| --- | --- | --- |
| $S_d$-N-30 | $f_{d,30}$ | 土日連休回数の下限 |
| $S_d$-N-31 | $f_{d,31}$ | 2連休回数の下限 |
| $S_d$-N-32 | $f_{d,32}$ | 3連続以上の休み回数の下限 |
| $S_d$-N-33 | $f_{d,33}$ | 単発日勤回数の上限 |
| $S_d$-N-34 | $f_{d,34}$ | ナース制約に関する違反度の最大最小化 |

## 4. 数値評価

急性期病院における病棟A，B，Cおよび慢性期病院における病棟Dにおいて，2024年11月分として実際に利用された制約条件および勤務希望表を用いて夜勤スケジューリングおよび日勤スケジューリングを行い，実行時間と最適値の観点から評価を行う．表6に対象病棟の概要を示す．ただし，看護師数における括弧内の数字は，夜勤可能な職員数を表す．また，利用した制約とは，病棟ごとに利用要否の選択が可能な「チーム/新人/男性/介護」の中から実際に利用した制約を表す．

表 6: 対象病棟の概要

| 病棟 | 看護師数 | 夜勤人数 | 配置シフト | 利用した制約 |
| --- | --- | --- | --- | --- |
| A | 29 (25) | 4（土日祝は3） | 日勤/夜勤/休み | 新人 |
| B | 39 (35) | 6 | 日勤/夜勤/休み | チーム，新人 |
| C | 91 (85) | 15 | 日勤/夜勤/休み | チーム，新人 |
| D | 30 (24) | 3 | 日勤/早出/遅出/夜勤/休み | チーム，男性，介護 |

まず，各病棟の夜勤スケジューリングおよび日勤スケジューリングにおける暫定解の推移を示す．今回は，Mac Studio (Apple M1 Max, 64 GB) において，Python 3.12.8 の Pulp 2.9.0 を用いて実行した．ソルバーとしては，CBC と CPLEX 22.1.0.0 を利用し，それらの実行ログを解析することで，実行時間と暫定解を取得した．

図 2 (a) から図 2 (f) に急性期病院である病棟 A，B，C に対する暫定解の推移を示す．病棟 A と比べて勤務人数や夜勤人数が多く，チーム制約を利用している病棟 B，C の方が，夜勤スケジューリングにおいて実行時間が長くなっている．特に，CBC を用いた場合にその差が顕著であり，最適解付近（$10^0 \sim 10^1$）から最適解に収束するまでの時間が長いことがわかる．一方で，日勤スケジューリングにおいては，実行時間における差は小さいものの，CBC の場合，夜勤スケジューリングと同様に最適解付近から最適解に収束するまでの時間が長くなっている．また，夜勤スケジューリングと日勤スケジューリングを比べると，夜勤スケジューリングの方が実行時間が長いことがわかる．

図 2 (g)，(h) に慢性期病院である病棟 D に対する暫定解の推移を示す．病棟 A，B，C とは異なり，夜勤スケジューリングよりも日勤スケジューリングの方が実行時間が長いことがわかる．これは，日勤スケジューリングで配置するシフトに早出と遅出を含むため，夜勤スケジューリングよりも解空間が広くなることに起因していると考えられる．

以上のように，病棟の規模や勤務形態などにより実行時間が異なるものの，いずれの病棟においても，実用的な時間[7]でスケジューリング問題の解が得られていることがわかる．

---

[7]後述する商用の自動作成支援サービスでは，運用上，1営業日以内に返答することが好ましいため，4病棟を契約している病院の場合，1〜2時間以内であれば実用的の範囲となる．

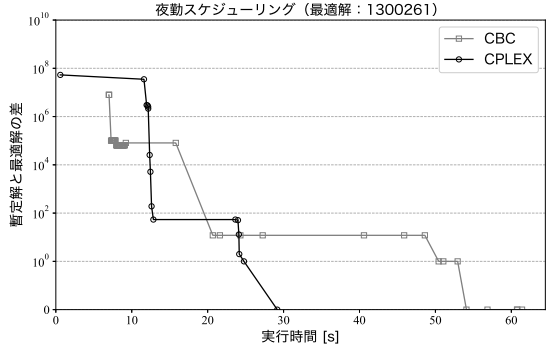

(a) 病棟 A：夜勤スケジューリング

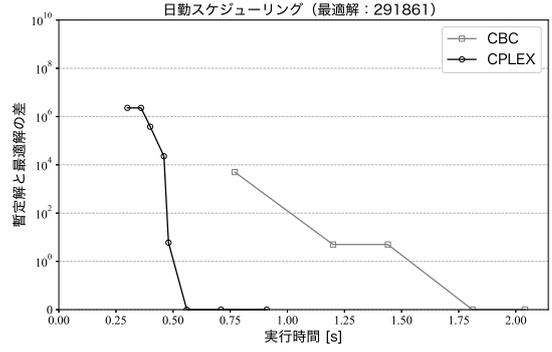

(b) 病棟 A：日勤スケジューリング

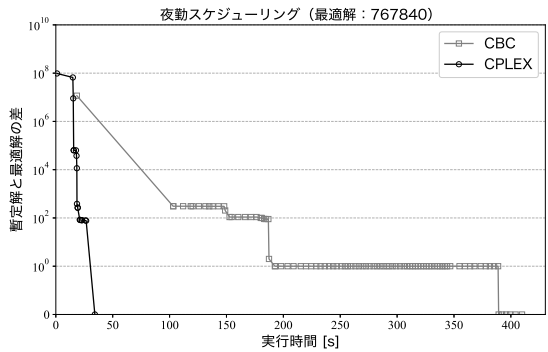

(c) 病棟 B 夜勤スケジューリング

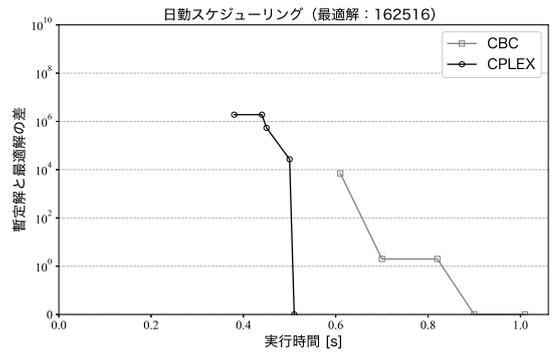

(d) 病棟 B 日勤スケジューリング

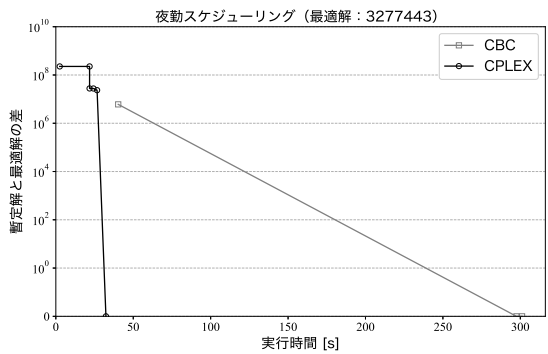

(e) 病棟 C：夜勤スケジューリング

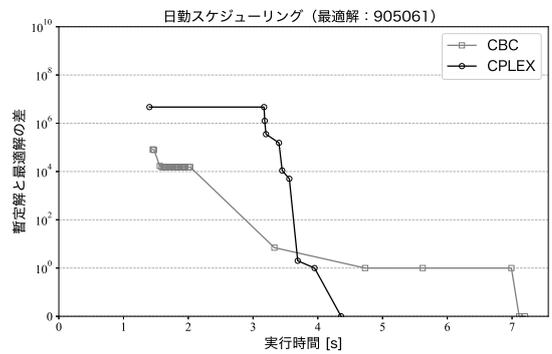

(f) 病棟 C：日勤スケジューリング

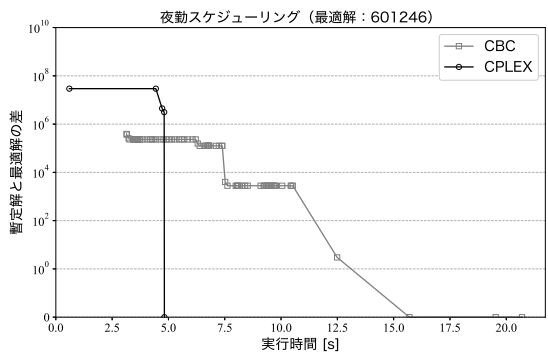

(g) 病棟 D：夜勤スケジューリング

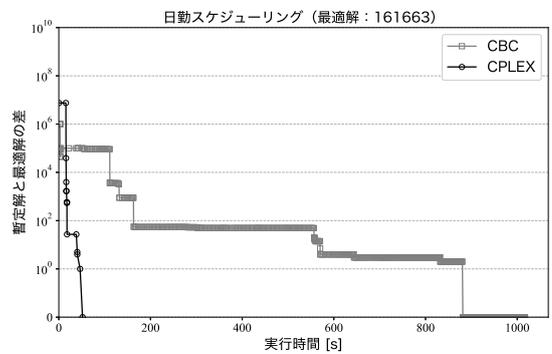

(h) 病棟 D：日勤スケジューリング

図 2: 各病棟における暫定解の推移

表 7: 看護師長の修正前後における最適値の変化

| 病棟 | 夜勤スケジューリング | | | 日勤スケジューリング | | |
| --- | --- | --- | --- | --- | --- | --- |
| | 修正前 | 修正後 | 変化率 | 修正前 | 修正後 | 変化率 |
| A | 1300261 | 1476259 | 1.14 | 291861 | 457789 | 1.57 |
| B | 767840 | 767840 | 1.00 | 162516 | 278791 | 1.72 |
| C | 3277443 | 4256446 | 1.30 | 905061 | 1189880 | 1.31 |
| D | 601246 | 917297 | 1.53 | 161663 | 616027 | 3.81 |

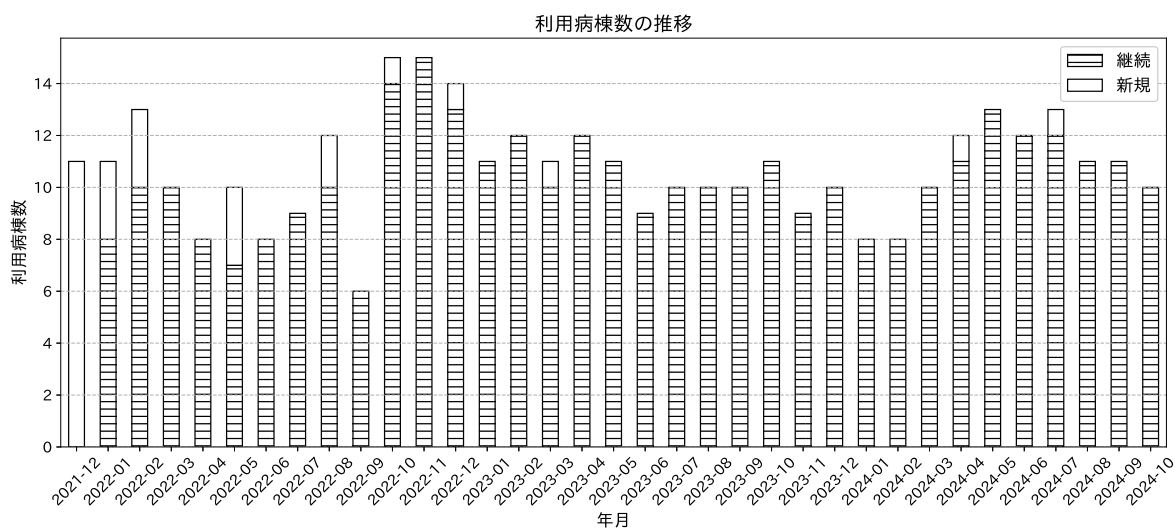

図 3: 大阪大学医学部附属病院における利用病棟数の推移

次に，看護師長が修正する前後における勤務表の最適値の変化を表 7 に示す．変化率（修正後/修正前）の大きさについては，目的関数の重みの設定に依存するためそれほど重要ではないが，少なくとも，看護師長が修正を行うことによって目的関数値が増加していることがわかる（病棟 B は夜勤スケジューリングの解に対する修正を行わなかったため変化なし）．したがって，看護師長にとって「納得感の高い」勤務表は，定式化された最適化問題において必ずしも最適解とはなっていないことがわかる．また，定式化が異なるため一概に比較することはできないが，すべての勤務を一括でスケジューリングする従来方法と比べて，2 段階スケジューリング法の方が，目的関数値自体は大きく（悪く）なると考えられる．

## 5. 事例紹介：急性期病院および慢性期病院における実運用

本節では，急性期病院である大阪大学医学部附属病院および慢性期病院である阪神リハビリテーション病院において，著者の仲島が 2 段階スケジューリング法を用いた自動作成支援を実際に運用した事例とそこから得られた課題について述べる．

### 5.1. 急性期病院の事例：大阪大学医学部附属病院

大阪大学医学部附属病院は，1086 床を有する急性期病院であり，2021 年 11 月（つまり，12 月の勤務表が対象）から試験運用を開始した．著者である仲島が所属しており，かつ，研究開発を目的としていることから，全病棟に対して一律に実施するのではなく，看護師長から要請のあった病棟に対してのみ自動作成支援を行ってきた．勤務表を作成する必要がある

病棟数は 26 であり，そのうち利用した病棟数の推移を図 3 に示す．平均して 10〜11 病棟程度（全体の 40% 程度）が利用していることがわかる．「新規」については，看護師長の異動や昇進に伴い，新たな病棟での利用があったことを示している．また，利用病棟数が変動しているのは，勤務表の作成条件が厳しい月のみ利用するといったスポット的な利用が発生するためである．ある程度の変動はあるが，全体としては安定して利用されており，特に，新しく昇格した看護師長を含む経験年数の浅い看護師長は継続的に利用している．

一方で，試験運用を開始した時点において，既に勤務表を 10 年以上作成している看護師長が半数程度を占めていた．そのようなベテラン看護師長は，勤務表作成自体に慣れており[8]，身体的・精神的負担が少ないため，積極的に自動作成支援を利用することはなかったと考えられる．このようなベテラン看護師長については，これまでの習慣を変えて自動作成支援を利用する方が負担になる可能性がある．また，これまでの慣習から勤務表作成は看護師長の責務であるという意識が強く，他者の介入を受け入れにくいという傾向も見受けられた．自動作成支援では，担当者のフィードバックコメント等も含むため，それを受け入れることが難しかった可能性も考えられる．この結果，最大利用数が 15 病棟となっているが，今後の世代交代により，継続的な利用者数が増加する可能性は高い．

### 5.2. 慢性期病院の事例：阪神リハビリテーション病院

阪神リハビリテーション病院は，192 床を有する慢性期病院であり，2023 年 6 月（つまり，7 月の勤務表が対象）から全 4 病棟において試験運用を開始した．そして，2024 年 9 月から，著者の仲島が代表を務める合同会社ナースラボ [12] の提供する有償サービス coNschel[9] (**cons**ulting for **N**urse **sche**duling) [13] としての利用に移行した．

試験運用にあたり，各看護師長から制約条件に関するヒアリングを行った．基本的には 3 節で示した大阪大学医学部附属病院と同様の制約条件を設定したが，ヒアリングの結果に合わせて以下の制約条件を変更した．

- 2 交代制 16 時間夜勤であるため，それに合わせて夜勤スケジューリングにおける夜勤インターバル等の複数の制約条件を調整した．
- 看護師とは異なる職能として介護職者が存在しており，毎日一人を夜勤に配置する必要があるため，新人看護師に関する勤務人数制約（$S_n$-S-17, $S_n$-S-18）と同形式の制約を追加した．
- 日勤スケジューリングにおいても同様に，介護職者に対して「日勤」，「早出」，「遅出」を毎日配置する必要があるため，新人看護師に関する勤務人数制約（$S_d$-S-17, $S_d$-S-18）と同形式の制約に加えて，「早出」および「遅出」に対する連続勤務に関する制約を追加した．

加えて，ヒアリング結果に合わせて目的関数の重みの設定を変更した．このように，慢性期病院の特性として介護職者を含むため，制約条件の設定が複雑化した．

運用当初においては，休み希望を締め切って看護師長が勤務表を作成している最中に追加の希望を依頼したり，勤務表作成後に不満な部分を看護師長に伝えて修正を依頼したりするといった，手作業で作成していた時の習慣が残っていたため，自動作成された勤務表を看護師長が大幅に修正した結果，勤務表作成に長時間を要するという問題が生じた．そこで，ある一人のスタッフの要望を優先することで，他のスタッフの $S_d$-N-26 から $S_d$-N-33 のナース制約違反が生じてしまう等の（数理最適化の観点から見た）勤務表の公平性や，休み希望

---

[8]中には，勤務表作成が好きな方もいた．
[9]コンシェルと読む．なお，Ns は看護師 (Nurse) の略記としてしばしば用いられる．

をハード制約とする場合，全スタッフのナース制約を全て満たす（または，同じだけ違反する）ことは現実的に不可能であるという制約の充足可能性に関連した限界に関する説明を全看護師長に対して行った．さらに，病棟間の公平性を高めるために，これまで病棟間で異なっていた，「希望の締切日，申請可能な希望シフトとその数，曜日ごとの配置人数，最大夜勤回数」等の運用ルールの標準化を行うように提案した．このような支援の結果，導入初期に生じた問題を解決することができ，勤務表作成にかかる時間が 70% 程度削減された．それに伴い，勤務希望を締め切ってから勤務表を公開するまでの期間が従来の 14〜20 日から約 7 日に短縮され，スタッフの満足度向上にも繋がった [20]．

また，阪神リハビリテーション病院では，看護部長が全病棟分の勤務表および勤務条件をまとめて担当者に共有するという形で運用している．このように看護部長が全病棟分を一括管理することで，全病棟の条件や現状を把握することができるため，病棟間の公平性を考慮した業務調整が可能となり，病棟間の応援体制（人数が不足している病棟に対し，一時的に人員を融通すること）に関する不公平感が軽減された．さらに，新規採用や退職があった場合の人事異動についても，これまでは各病棟の看護師長の声に依存することが多かったが，運用ルールを標準化し，公平化を促進したことで，客観的な判断が可能となった [16, 20]．

### 5.3. 実運用から得られた課題と対策

特徴の異なる二つの病院における実運用を通じて，以下の三つの課題が明らかになった．

A. 勤務表作成自体に負担を感じていない看護師長に対しては，担当者とのやりとり（自動作成支援）が逆に手間となってしまうため，積極的な利用を促すことが難しい．

B. 手作業時の運用から脱却できない場合，導入効果が薄くなる．

C. 担当者を介する必要があるため，その数がスケールアップのボトルネックとなる．

課題 A については，人事異動や応援体制等の病棟間の公平性を向上させるという観点から，全病棟に対してトップダウンで導入するという方法が考えられるが，ベテラン看護師長からの反発が容易に予想される．そこで今後は，生成 AI を用いた自動作成支援用の音声対話可能なチャットボットを開発し，担当者とのやりとりを自動化することを検討している（自動で条件設定まで行うことが目標）．さらに，多様な解の生成に関する研究 [6] を参考にして，看護師長が直感的に好みの勤務表を選択できるような UI を開発することも検討している．

また，自動作成支援サービスを提供する企業については，サービスの導入および展開に際して，課題 B ならびに C に対する対策を講じることが重要である．課題 B については，システムによって解決できる問題ではなく，各病棟の文化を変えることが必要となる．そのため，阪神リハビリテーション病院のような事例を一つでも多く作り，自動作成支援の質を高めていくことが重要である．課題 C については，上記の AI チャットボット等の開発に加えて，数理最適化および看護の両方の知識を持つ人材を育成することが重要である．現場の看護師の気持ちを理解できなければ，円滑なコミュニケーションおよび自動作成支援の継続は難しいと考えられる．

### 6. おわりに

本稿では，2 段階スケジューリング法を用いた看護師勤務表の自動作成支援について，夜勤スケジューリングおよび日勤スケジューリングの定式化を行うとともに，急性期病院および慢性期病院という特徴の異なる二つの病院において実運用した結果について述べた．その実運用の結果，自動作成支援によって看護師長の負担が軽減されるだけではなく，病棟間の応援体制や人事異動等の管理面にも効果があることが示された．2 段階スケジューリング法

は，看護現場に限らず，医療や介護に加え，救急や消防，警備，インフラ保守など，24 時間の人員配置が必要となり，かつ，様々な暗黙知を前提としたスケジューリングが行われている他業種にも応用可能であると考えられる．一方で，これまでの手作業の文化を変えることが難しいという課題やスケールアップが難しいという課題が明らかになった．合同会社ナースラボだけではなく，各地域において同様のサービスを提供する企業が増えることで，より勤務表の自動作成が普及することが期待される．

生成 AI の登場によってプログラミングのハードルが低くなった今では，各病院において自前で勤務表自動作成システムを開発・運用することも可能である．医療情報部門や事務部門における技術者と連携し，本稿における定式化の詳細と数理最適化およびプログラミングに関する書籍等 [4, 5] を参考にして，その実現を目指す看護部の手助けとなれば幸いである．

**参考文献**

仲島圭将
大阪大学 大学院医学系研究科 医療情報学
〒 565-0871 大阪府吹田市山田丘 2-2
E-mail: k-nakashima@hp-nurse.med.osaka-u.ac.jp

## A. 夜勤スケジューリングにおけるハード制約の詳細

**$H_n$-N-1:** 勤務希望と指定勤務の固定および夜勤シフトのみに限定

看護師 $n$ の日付 $d$ に対し，本人が希望する勤務シフトまたは会議や研修等の理由で指定された勤務シフト $s$ が与えられた場合（$x'_{n,d,s}=1$），そのシフトを固定する（$x_{n,d,s}=1$）．また，夜勤スケジューリングでは「入，明，休，未」のみを配置するため，それら以外のシフトは配置しないように固定する（$x_{n,d,s}=0$）．

$$\begin{cases} x_{n,d,s} = 1 & \text{if } x'_{n,d,s} = 1 \\ x_{n,d,s} = 0 & \text{if } x'_{n,d,s} = 0,\ s \notin \{\text{入},\text{明},\text{休},\text{未}\} \end{cases} \quad (n \in \mathcal{N},\ d \in \mathcal{D},\ s \in \mathcal{S})$$

**$H_n$-N-2:** 前月からの月跨ぎの禁止シフト

前月末から月跨ぎで続く 3 連続の日勤シフト（$s \in \mathcal{S}_{\text{day}}$）に対して，「未–入–明」が連続して配置された場合，「未」が夜勤スケジューリングの後工程において「12h」に置換されて 6 連勤となるため禁止する．

$$\sum_{h=0}^{2} \sum_{s \in \mathcal{S}_{\text{day}}} x_{n,d+h,s} + x_{n,d+3,\text{未}} + x_{n,d+4,\text{入}} \leq 4 \quad (n \in \mathcal{N}_{\text{n}},\ d \in \mathcal{D}_{\text{p}})$$

**$H_n$-N-3:** 前月からの月跨ぎの夜勤シフト固定

前月末の勤務が「入」の場合，対象月の 1 日は「明」に固定する．前月末の勤務が「入」ではない場合，対象月の 1 日に「明」が割り当てられることはない．

$$x_{n,d_1,\text{明}} = x'_{n,d_{\text{p,last}},\text{入}} \quad (n \in \mathcal{N}_{\text{n}} \cup \mathcal{N}_{\text{no}})$$

**$H_n$-N-4:** 次月への月跨ぎの固定・禁止シフト

次月の 1 日に休み希望がある場合，対象月の最終日に「入」を割り当てない．また，次月の 1 日に「明」が割り当てられている場合，対象月の最終日を「入」に固定する．

$$\begin{cases} x_{n,d_{\text{last}},\text{入}} = 0 & \text{if } x'_{n,d_{\text{n1}},\text{休}} + x'_{n,d_{\text{n1}},\text{特休}} \geq 1 \\ x_{n,d_{\text{last}},\text{入}} = 1 & \text{if } x'_{n,d_{\text{n1}},\text{明}} = 1 \end{cases} \quad (n \in \mathcal{N}_{\text{n}} \cup \mathcal{N}_{\text{no}})$$

12 時間夜勤の場合,「12h–入–明–休」が基本となるため，次月の 1, 2 日に日勤や会議等がある場合，対象月の最終日に「入」および「明」を割り当てない．次月の 2 日に「入」が割り当てられている場合，対象月の最終日に「入」および「明」を割り当てない．また，次月の 2 日に「明」が割り当てられている場合，次月の 1 日を「入」に固定し，対象月の最終日に「入」および「明」を割り当てない．

$$\begin{cases} x_{n,d_{\text{last}},\text{入}} = 0,\ x_{n,d_{\text{last}},\text{明}} = 0 & \text{if } x'_{n,d_{\text{n1}},\text{日}} + x'_{n,d_{\text{n1}},\text{他}} \geq 1 \\ x_{n,d_{\text{last}},\text{入}} = 0 & \text{if } x'_{n,d_{\text{n2}},\text{日}} + x'_{n,d_{\text{n2}},\text{他}} \geq 1 \\ x_{n,d_{\text{last}},\text{入}} = 0,\ x_{n,d_{\text{last}},\text{明}} = 0 & \text{if } x'_{n,d_{\text{n2}},\text{入}} = 1 \\ x_{n,d_{\text{last}},\text{入}} = 0,\ x_{n,d_{\text{last}},\text{明}} = 0,\ x_{n,d_{\text{n1}},\text{入}} = 1 & \text{if } x'_{n,d_{\text{n2}},\text{明}} = 1 \end{cases} \quad (n \in \mathcal{N}_{\text{n}})$$

**H$_\text{n}$-N-5:** 6 連勤の禁止

前月末からの月跨ぎも含め，連続する 6 日間に割り当てられる勤務シフトの数を 5 以下にする．

$$\sum_{h=0}^{5} \sum_{s \in \mathcal{S}_\text{work}} x_{n,d+h,s} \leq 5 \quad (n \in \mathcal{N},\ d \in \mathcal{D}_\text{p} \cup \mathcal{D}_\text{t})$$

**H$_\text{n}$-N-6:** 禁止シフト

「入–休」や「入–日」など，禁止シフト集合に含まれるシフト系列を割り当てない．

$$\sum_{i=0}^{|\boldsymbol{s}^{(\text{h})}|-1} x_{n,d+i,s_i^{(\text{h})}} \leq |\boldsymbol{s}^{(\text{h})}| - 1 \quad (n \in \mathcal{N},\ d \in \mathcal{D}_\text{f},\ \boldsymbol{s}^{(\text{h})} \in \mathcal{S}_\text{hf})$$

ここで，$s_i^{(\text{h})}$ は禁止シフト系列 $\boldsymbol{s}^{(\text{h})}$ の $i$ 番目の要素を表す．例えば，$s_i^{(\text{h})} = \{\,入, 休\,\}$ の時，$s_1^{(\text{h})} = 入$，$s_2^{(\text{h})} = 休$ となる．実際には，看護師の勤務形態（夜勤専従，日勤専従など）に応じて禁止シフト集合を設定する必要があることに注意する．

**H$_\text{n}$-N-7:** 日勤専従看護師の夜勤制限

日勤専従看護師 $n\,(\in \mathcal{N}_\text{day})$ に「入，明」を割り当てない．

$$x_{n,d,\,入} = 0,\ x_{n,d,\,明} = 0 \quad (n \in \mathcal{N}_\text{day},\ d \in \mathcal{D}_\text{t})$$

**H$_\text{n}$-N-8:** 夜勤可能および夜勤専従看護師の夜勤回数制限

夜勤に従事する看護師 $n\,(\in \mathcal{N}_\text{n} \cup \mathcal{N}_\text{no})$ の夜勤回数を制限する必要がある場合（$f_n^{(\text{n})} = 1$），指定した回数の範囲に収まるようにする．

$$\begin{cases} N_n^{(\text{LB})} \leq \displaystyle\sum_{d \in \mathcal{D}} x_{n,d,\,入} \leq N_n^{(\text{UB})} \\ N_n^{(\text{LB})} \leq \displaystyle\sum_{d \in \mathcal{D}} x_{n,d,\,明} \leq N_n^{(\text{UB})} \end{cases} \text{if } f_n^{(\text{n})} = 1 \qquad (n \in \mathcal{N}_\text{n} \cup \mathcal{N}_\text{no})$$

**H$_\text{n}$-N-9:** 夜勤インターバルの下限

12 時間夜勤の場合，「12h–入–明–休」が基本となるため，「入」に関して最低 4 日以上のインターバルを確保する．つまり，連続する 4 日間に割り当てられる「入」の回数を 1 以下にする．

$$\sum_{h=0}^{3} x_{n,d+h,\,入} \leq 1 \quad (n \in \mathcal{N}_\text{n},\ d \in \mathcal{D}_\text{ilb})$$

一方，16 時間夜勤の場合，この制約は適用しない．

**H$_\text{n}$-N-10:** 同一勤務の連続回数制限

勤務シフト $s\,(\in \mathcal{S}_\text{work})$ に対して，連続回数が最大回数 $N_{s,\text{max}}$ 以下になるようにする．例えば，$N_{日,\text{max}} = 5$, $N_{入,\text{max}} = 1$ となる．

$$\sum_{h=0}^{N_{s,\text{max}}} x_{n,d+h,s} \leq N_{s,\text{max}} \quad (n \in \mathcal{N},\ d \in \mathcal{D}_s,\ s \in \mathcal{S}_\text{work})$$

**H$_\text{n}$-N-11:** 3 連休候補を 2 箇所以上確保

休みに関する好みが「3 連休」の看護師（$f_{n,3\,連休} = 1$）に対し，日勤スケジューリング時に 3 連休を割り当てることができるように，2 箇所以上の候補を確保する．つまり，「休，特休，未」が 3 連続で続いた後に，「入，明」以外のシフトが割り当てられている場合を 3 連休候補とし，日付 $d$ を初日とする 3 連休候補が成立するか否かを表す二値変数 $s_{n,d}^{(3\text{off})} \in \{0,1\}$ の合計が 2 以上になるようにする．こ

れは，12 時間夜勤の場合，後工程において「入」の前日の「未」が「12h」に置換されるためである．

$$\begin{cases} \sum_{h=0}^{2}\sum_{s\in\mathcal{S}'_{\text{off}}} x_{n,d+h,s} + \sum_{s\in\mathcal{S}_{\text{day}}\cup\mathcal{S}'_{\text{off}}} x_{n,d+3,s} \geq 4 \cdot s_{n,d}^{(3\text{off})} \\ \sum_{h=0}^{2}\sum_{s\in\mathcal{S}'_{\text{off}}} x_{n,d+h,s} + \sum_{s\in\mathcal{S}_{\text{day}}\cup\mathcal{S}'_{\text{off}}} x_{n,d+3,s} \leq 3 + s_{n,d}^{(3\text{off})} \end{cases} \quad (n\in\mathcal{N},\ d\in\mathcal{D}_{3\text{off}})$$

$$\begin{cases} \sum_{d\in\mathcal{D}_{3\text{off}}} s_{n,d}^{(3\text{off})} \geq 2 & \text{if } f_{n,3\text{ 連休}} = 1 \\ \sum_{d\in\mathcal{D}_{3\text{off}}} s_{n,d}^{(3\text{off})} \geq 0 & \text{otherwise} \end{cases} \quad (n\in\mathcal{N})$$

**H$_\text{n}$-S-12:** 日勤においてリーダー業務可能な看護師を一人以上確保

日勤スケジューリングにおいて，リーダー業務可能な看護師を一人以上確保できるようにするため，「日，未」のいずれかのシフトが割り当てられたリーダー可能な看護師を一人以上確保する．ただし，土日祝においては勤務人数が少ないため，12 時間日勤に従事する看護師がリーダー業務を担当することを考慮する必要がある．

$$\begin{cases} \sum_{n\in\mathcal{N}_\text{L}}\sum_{s\in\{\text{日, 未}\}} x_{n,d,s} \geq 1 & \text{if weekday}(d) = 1,\ \text{i.e., 平日} \\ \sum_{n\in\mathcal{N}_\text{L}}\sum_{s\in\{\text{日,12h, 未}\}} x_{n,d,s} \geq 1 & \text{if weekday}(d) = 0,\ \text{i.e., 土日祝} \end{cases} \quad d\in\mathcal{D}_\text{t}$$

**H$_\text{n}$-S-13:** 男性看護師の人数制限

病棟の特性によっては，必ず一人以上の男性看護師を配置する必要がある場合や，女性患者への配慮から男性看護師の配置上限を設ける場合がある．

$$N_{\text{n,men}}^{(\text{LB})} \leq \sum_{n\in\mathcal{N}_\text{men}} x_{n,d,\text{入}} \leq N_{\text{n,men}}^{(\text{UB})} \quad d\in\mathcal{D}_\text{t}$$

**H$_\text{n}$-S-14:** 必ず勤務記号が割り当てられる

必ず一つの勤務記号（「未」を含む）が割り当てられるようにする．

$$\sum_{s\in\mathcal{S}} x_{n,d,s} = 1 \quad (n\in\mathcal{N},\ d\in\mathcal{D})$$

## B. 夜勤スケジューリングにおけるソフト制約の詳細

**S$_\text{n}$-S-1:** 日毎の勤務人数の下限

看護師数が下限を下回る場合，下限との差をペナルティとして与える．

$$f_{\text{n},1}(\boldsymbol{X}) = \sum_{d\in\mathcal{D}_\text{t}} \max\big(0, N_\text{n}^{(\text{LB})} - \sum_{n\in\mathcal{N}_\text{grp}} x_{n,d,\text{入}}\big)$$

**S$_\text{n}$-S-2:** 日毎の勤務人数の上限

看護師数が上限を超える場合，上限との差をペナルティとして与える．

$$f_{\text{n},2}(\boldsymbol{X}) = \sum_{d\in\mathcal{D}_\text{t}} \max\big(0, \sum_{n\in\mathcal{N}_\text{grp}} x_{n,d,\text{入}} - N_\text{n}^{(\text{UB})}\big)$$

以下の $f_{\text{n},3}$ から $f_{\text{n},16}$ は，$f_{\text{n},1}$，$f_{\text{n},2}$ のグループおよびチーム拡張であるため，その形式は同様である．また，$f_{\text{n},17}$，$f_{\text{n},18}$ は，新人看護師を対象とした制約条件ではあるが，これらについてもその形式は同様である．

**S$_\text{n}$-S-3, S$_\text{n}$-S-4:** グループ 1 に属する看護師の勤務人数の下限・上限

$$f_{\text{n},3}(\boldsymbol{X}) = \sum_{d\in\mathcal{D}_\text{t}} \max\big(0, N_{\text{n},1}^{(\text{LB})} - \sum_{n\in\mathcal{N}_1} x_{n,d,\text{入}}\big) \quad f_{\text{n},4}(\boldsymbol{X}) = \sum_{d\in\mathcal{D}_\text{t}} \max\big(0, \sum_{n\in\mathcal{N}_1} x_{n,d,\text{入}} - N_{\text{n},1}^{(\text{UB})}\big)$$

**S$_\text{n}$-S-5, S$_\text{n}$-S-6:** グループ 1，2 に属する看護師の勤務人数の下限・上限

$$f_{\text{n},5}(\boldsymbol{X}) = \sum_{d\in\mathcal{D}_\text{t}} \max\big(0, N_{\text{n},12}^{(\text{LB})} - \sum_{n\in\mathcal{N}_{12}} x_{n,d,\text{入}}\big) \quad f_{\text{n},6}(\boldsymbol{X}) = \sum_{d\in\mathcal{D}_\text{t}} \max\big(0, \sum_{n\in\mathcal{N}_{12}} x_{n,d,\text{入}} - N_{\text{n},12}^{(\text{UB})}\big)$$

**S$_\text{n}$-S-7, S$_\text{n}$-S-8:** グループ 4 に属する看護師の勤務人数の下限・上限

$$f_{\text{n},7}(\boldsymbol{X}) = \sum_{d \in \mathcal{D}_\text{t}} \max\bigl(0, N_{\text{n},4}^{(\text{LB})} - \sum_{n \in \mathcal{N}_4} x_{n,d,\text{入}}\bigr) \quad f_{\text{n},8}(\boldsymbol{X}) = \sum_{d \in \mathcal{D}_\text{t}} \max\bigl(0, \sum_{n \in \mathcal{N}_4} x_{n,d,\text{入}} - N_{\text{n},4}^{(\text{UB})}\bigr)$$

**S$_\text{n}$-S-9, S$_\text{n}$-S-10:** チーム $t$ およびグループ 1 に属する看護師の勤務人数の下限・上限

$$f_{\text{n},9}(\boldsymbol{X}) = \sum_{t \in \mathcal{T}} \sum_{d \in \mathcal{D}_\text{t}} \max\bigl(0, N_{\text{n},t,1}^{(\text{LB})} - \sum_{n \in \mathcal{N}_{t,1}} x_{n,d,\text{入}}\bigr)$$

$$f_{\text{n},10}(\boldsymbol{X}) = \sum_{t \in \mathcal{T}} \sum_{d \in \mathcal{D}_\text{t}} \max\bigl(0, \sum_{n \in \mathcal{N}_{t,1}} x_{n,d,\text{入}} - N_{\text{n},t,1}^{(\text{UB})}\bigr)$$

**S$_\text{n}$-S-11, S$_\text{n}$-S-12:** チーム $t$ かつグループ 1, 2 に属する看護師の勤務人数の下限・上限

$$f_{\text{n},11}(\boldsymbol{X}) = \sum_{t \in \mathcal{T}} \sum_{d \in \mathcal{D}_\text{t}} \max\bigl(0, N_{\text{n},t,12}^{(\text{LB})} - \sum_{n \in \mathcal{N}_{t,12}} x_{n,d,\text{入}}\bigr)$$

$$f_{\text{n},12}(\boldsymbol{X}) = \sum_{t \in \mathcal{T}} \sum_{d \in \mathcal{D}_\text{t}} \max\bigl(0, \sum_{n \in \mathcal{N}_{t,12}} x_{n,d,\text{入}} - N_{\text{n},t,12}^{(\text{UB})}\bigr)$$

**S$_\text{n}$-S-13, S$_\text{n}$-S-14:** チーム $t$ かつグループ 3, 4 に属する看護師の勤務人数の下限・上限

$$f_{\text{n},13}(\boldsymbol{X}) = \sum_{t \in \mathcal{T}} \sum_{d \in \mathcal{D}_\text{t}} \max\bigl(0, N_{\text{n},t,34}^{(\text{LB})} - \sum_{n \in \mathcal{N}_{t,34}} x_{n,d,\text{入}}\bigr)$$

$$f_{\text{n},14}(\boldsymbol{X}) = \sum_{t \in \mathcal{T}} \sum_{d \in \mathcal{D}_\text{t}} \max\bigl(0, \sum_{n \in \mathcal{N}_{t,34}} x_{n,d,\text{入}} - N_{\text{n},t,34}^{(\text{UB})}\bigr)$$

**S$_\text{n}$-S-15, S$_\text{n}$-S-16:** チーム $t$ に属する看護師の勤務人数の下限

$$f_{\text{n},15}(\boldsymbol{X}) = \sum_{t \in \mathcal{T}} \sum_{d \in \mathcal{D}_\text{t}} \max\bigl(0, N_{\text{n},t}^{(\text{LB})} - \sum_{n \in \mathcal{N}_t} x_{n,d,\text{入}}\bigr)$$

$$f_{\text{n},16}(\boldsymbol{X}) = \sum_{t \in \mathcal{T}} \sum_{d \in \mathcal{D}_\text{t}} \max\bigl(0, \sum_{n \in \mathcal{N}_t} x_{n,d,\text{入}} - N_{\text{n},t}^{(\text{UB})}\bigr)$$

**S$_\text{n}$-S-17, S$_\text{n}$-S-18:** 新人看護師の勤務人数の下限・上限

$$f_{\text{n},17}(\boldsymbol{X}) = \sum_{d \in \mathcal{D}_\text{t}} \max\bigl(0, N_{\text{r}}^{(\text{LB})} - \sum_{n \in \mathcal{N}_\text{r}} x_{n,d,\text{入}}\bigr) \quad f_{\text{n},18}(\boldsymbol{X}) = \sum_{d \in \mathcal{D}_\text{t}} \max\bigl(0, \sum_{n \in \mathcal{N}_\text{r}} x_{n,d,\text{入}} - N_{\text{r}}^{(\text{UB})}\bigr)$$

**S$_\text{n}$-S-19:** ペア勤務回数の下限

プリセプターとプリセプティのような，業務上，月に数回程度は同日に勤務させたい関係をペア関係という．ペア関係 $p_i$ にある看護師が日付 $d$ に指定した勤務シフトに割り当てられたか否かを表す二値変数 $s_{i,d}^{(\text{pair})}$ の合計が下限 $N_{i,\min}$ 未満である場合，下限との差をペナルティとして与える．

$$f_{\text{n},19}(\boldsymbol{X}) = \sum_{p_i \in \mathcal{P}} \max\bigl(0, N_{i,\min} - \sum_{d \in \mathcal{D}_\text{t}} s_{i,d}^{(\text{pair})}\bigr)$$

ただし，$p_i = (n_1, n_2, s_1, s_2, N_{i,\min}) \in \mathcal{P}$, $d \in \mathcal{D}_\text{t}$ に対し，以下を満たす．

$$x_{n_1,d,s_1} + x_{n_2,d,s_2} \geq 2 \cdot s_{i,d}^{(\text{pair})}, \quad x_{n_1,d,s_1} + x_{n_2,d,s_2} \leq s_{i,d}^{(\text{pair})} + 1$$

**S$_\text{n}$-S-20:** 禁止ペア

能力や性格などの理由から同日に勤務させたくない関係を禁止ペア関係という．禁止ペア関係 $p_i^{(\text{f})} = (n_1, n_2, s_1, s_2) \in \mathcal{P}_\text{f}$ にある看護師が同日 $d$ に指定した勤務シフトに割り当てられた場合にペナルティを与える．

$$f_{\text{n},20}(\boldsymbol{X}) = \sum_{p_i^{(\text{f})} \in \mathcal{P}_\text{f}} \sum_{d \in \mathcal{D}_\text{t}} \max\bigl(0, (x_{n_1,d,s_1} + x_{n_2,d,s_2} - 1)\bigr)$$

**S$_n$-N-21:** 禁止勤務
家庭の事情等により，特定の曜日に対して勤務制限がある場合，その曜日に対して指定したシフトが割り当てられた場合にペナルティを与える，つまり，禁止勤務集合 $\mathcal{F}$ に含まれる勤務シフト $f_i = (n, s, d_\mathrm{w})$ が割り当てられた場合にペナルティを与える．

$$f_{\mathrm{n},21}(\boldsymbol{X}) = \sum_{f_i \in \mathcal{F}} \sum_{\substack{d \in \mathcal{D}_\mathrm{t}, \\ \mathrm{dayW}(d) = d_\mathrm{w}}} x_{n,d,s}$$

ここで，$\mathrm{dayW}(d)$ は日付 $d$ の曜日（祝日含む）を表す．

**S$_n$-N-22:** ソフト禁止シフト
ハード制約条件における禁止シフト制約（H$_n$-N-6）を緩和した制約であり，基本的には配置すべきではないが，状況によっては配置せざるを得ないシフト系列 $\boldsymbol{s}_i^{(s)}$ が割り当てられた場合にペナルティを与える．

$$f_{\mathrm{n},22}(\boldsymbol{X}) = \sum_{n \in \mathcal{N}} \sum_{\boldsymbol{s}^{(s)} \in \mathcal{S}_\mathrm{sf}} \sum_{d \in \mathcal{D}_\mathrm{f}} \max\left(0, \sum_{i=0}^{|\boldsymbol{s}^{(s)}|-1} x_{n,d+i,s_i^{(s)}} - (|\boldsymbol{s}^{(s)}| - 1)\right)$$

ここで，$s_i^{(\mathrm{s})}$ は禁止シフト系列 $\boldsymbol{s}^{(\mathrm{s})}$ の $i$ 番目の要素を表す．

**S$_n$-N-23:** 休み（週休，祝日，振替）回数を 4 回以上残す
夜勤（12h–入–明–休）が多く割り当てられると，その分だけ日勤スケジューリング時に配置可能な休み回数が減少してしまい，休みの好みに合わせたスケジューリングが困難になる．そこで，配置可能な休み回数を 4 回以上残すようにするため，対象月に配置可能な休み回数 $N_\mathrm{off}$ と配置された休み回数の差が 4 未満の場合，その分だけペナルティとして与える．

$$f_{\mathrm{n},23}(\boldsymbol{X}) = \sum_{n \in \mathcal{N}} \max\left(0, 4 - (N_\mathrm{off} - \sum_{d \in \mathcal{D}_\mathrm{t}} x_{n,d,\text{休}})\right)$$

**S$_n$-N-24:** 夜勤インターバルの上限
サーカディアンリズム [10] の観点から，夜勤の配置間隔をなるべく一定に保つことが望ましい．そのため，夜勤インターバルが上限の 14 日を超える，つまり，連続する14日間に一度も「入」が割り当てられていない場合にペナルティを与える．

$$f_{\mathrm{n},24}(\boldsymbol{X}) = \sum_{n \in \mathcal{N}_\mathrm{n} \cup \mathcal{N}_\mathrm{no}} \sum_{d \in \mathcal{D}_\mathrm{iub}} \max\left(0, 1 - \sum_{h=0}^{13} x_{n,d+h,\text{入}}\right)$$

**S$_n$-N-25:** 夜勤回数の下限
「入，明」は必ず続けて配置されるが，月を跨いで配置された場合，対象月の「入」と「明」の合計回数が一致しないことがある．よって，「入」と「明」の合計回数のうち，小さい方が下限回数 $N_n^{(\mathrm{LB})}$ 未満の場合，その差をペナルティとして与える．

$$f_{\mathrm{n},25}(\boldsymbol{X}) = \sum_{n \in \mathcal{N}_\mathrm{n} \cup \mathcal{N}_\mathrm{no}} \max\left(0, N_n^{(\mathrm{LB})} - \min\left(\sum_{d \in \mathcal{D}_\mathrm{t}} x_{n,d,\text{入}}, \sum_{d \in \mathcal{D}_\mathrm{t}} x_{n,d,\text{明}}\right)\right)$$

**S$_n$-N-26:** 夜勤回数の上限
$f_{\mathrm{n},25}$ と同様に，「入」と「明」の合計回数のうち，大きい方が上限回数 $N_n^{(\mathrm{UB})}$ を超える場合，その差をペナルティとして与える．

$$f_{\mathrm{n},26}(\boldsymbol{X}) = \sum_{n \in \mathcal{N}_\mathrm{n} \cup \mathcal{N}_\mathrm{no}} \max\left(0, \max\left(\sum_{d \in \mathcal{D}_\mathrm{t}} x_{n,d,\text{入}}, \sum_{d \in \mathcal{D}_\mathrm{t}} x_{n,d,\text{明}}\right) - N_n^{(\mathrm{UB})}\right)$$

**S$_n$-N-27:** 金土日月の夜勤回数の公平化
日勤スケジューリングにおいて，土日休みを公平に割り当てることができるように，金土日月の夜勤回数を平滑化する．12時間夜勤の場合，「12h–入–明–休」が基本となっているため，金曜日および月曜日が含まれていることに注意する．金土日月における合計回数と平均回数 $\bar{N}_\mathrm{fssm}$ との差をペナルティとして与える．このとき，休みの好みが「土日連休」の看護師に対してはペナルティの重みを大きくする．

$$f_{\mathrm{n},27}(\boldsymbol{X}) = \sum_{n \in \mathcal{N}_\mathrm{n}} w_{27,n}^{(\mathrm{n})} \max\left(0, \sum_{d \in \mathcal{D}_\mathrm{fssm}} x_{n,d,\text{入}} - \bar{N}_\mathrm{fssm}\right)$$

$$w_{27,n}^{(\mathrm{n})} = \begin{cases} 5 & \text{if } f_n^{(\mathrm{p})} = \text{土日連休} \\ 1 & \text{otherwise} \end{cases}$$

**S$_\text{n}$-N-28:** 手術日等の夜勤回数の公平化

急性期病院においては，手術や検査等のイベントの曜日が固定されていることが多く，その曜日の夜勤業務は負担が大きい．そこで，その曜日の夜勤回数を平滑化するため，該当曜日の合計回数と平均回数 $\bar{N}_\text{ope}$ との差をペナルティとして与える．

$$f_{\text{n},28}(\boldsymbol{X}) = \sum_{n \in \mathcal{N}_\text{n}} \max\bigl(0, \sum_{d \in \mathcal{D}_\text{ope}} x_{n,d,\,\text{入}} - \bar{N}_\text{ope}\bigr)$$

**S$_\text{n}$-N-29:** 夜勤を含む 5 連続勤務の回避

3 日間連続して日勤等の勤務が続いた後に「入–明」が割り当てられるという非常に負担の大きい 5 連続勤務に対してペナルティを与える．12 時間夜勤の場合，「入」の前日の「未」が後工程において「12h」に置換されるため，3 日目に「未」が配置された場合も考慮する．また，休みの好みが「単発休み」，つまり，連続勤務を好まない看護師に対してはペナルティの重みを大きくする．

$$f_{\text{n},29}(\boldsymbol{X}) = \sum_{n \in \mathcal{N}_\text{n}} \sum_{d \in \mathcal{D}_\text{5d}} w^{(\text{n})}_{29,n} \max\bigl(0, \sum_{h=0}^{1} \sum_{s \in \mathcal{S}_\text{day}} x_{n,d+h,s} + \sum_{s \in \{\text{12h, 未}\}} x_{n,d+2,s} + x_{n,d+3,\,\text{入}} - 3\bigr)$$

$$w^{(\text{n})}_{29,n} = \begin{cases} 5 & \text{if } f^{(\text{p})}_n = \text{単発休み} \\ 1 & \text{otherwise} \end{cases}$$

**S$_\text{n}$-N-30:** 夜勤を含む「5 連勤–1 休み–4 連勤」および「4 連勤–1 休み–5 連勤」の回避

$f_{\text{n},29}$ と同様に，夜勤を含む 4 連勤または 5 連勤が単発休みを挟んで続くことは非常に負担が大きいため，いずれかが割り当てられた場合にペナルティを与える．

$$f_{\text{n},30}(\boldsymbol{X}) = \sum_{n \in \mathcal{N}_\text{n}} \sum_{d \in \mathcal{D}_\text{5d}} \max\bigl(0, \max(f^{(\text{A})}_{\text{n},30,n}, f^{(\text{B})}_{\text{n},30,n}) - 9\bigr)$$

$$f^{(\text{A})}_{\text{n},30,n} = \sum_{h=0}^{1} \sum_{s \in \mathcal{S}_\text{day}} x_{n,d+h,s} + \sum_{s \in \{\text{12h, 未}\}} x_{n,d+2,s} + x_{n,d+3,\,\text{入}} + x_{n,d+4,\,\text{明}}$$
$$+ \sum_{s \in \mathcal{S}_\text{off}} x_{n,d+5,s} + \sum_{s \in \mathcal{S}_\text{day}} x_{n,d+6,s} + \sum_{s \in \{\text{12h, 未}\}} x_{n,d+7,s} + x_{n,d+8,\,\text{入}} + x_{n,d+9,\,\text{明}}$$

$$f^{(\text{B})}_{\text{n},30,n} = \sum_{s \in \mathcal{S}_\text{day}} x_{n,d,s} + \sum_{s \in \{\text{12h, 未}\}} x_{n,d+1,s} + x_{n,d+2,\,\text{入}} + x_{n,d+3,\,\text{明}} + \sum_{s \in \mathcal{S}_\text{off}} x_{n,d+4,s}$$
$$+ \sum_{h=5}^{6} \sum_{s \in \mathcal{S}_\text{day}} x_{n,d+h,s} + \sum_{s \in \{\text{12h, 未}\}} x_{n,d+7,s} + x_{n,d+8,\,\text{入}} + x_{n,d+9,\,\text{明}}$$

以下の $f_{\text{n},31}$ から $f_{\text{n},33}$ は，夜勤スケジューリングにおいて，夜勤が連続して割り当てられることを避ける，つまり，一定間隔を空けて配置できるようにするための制約である．

**S$_\text{n}$-N-31:** 3 連続夜勤の回避

12 時間夜勤（12h–入–明–休）に従事する夜勤可能看護師に対し，連続する 9 日間で 3 回の夜勤が割り当てられた場合にペナルティを与える．このとき，休みの好みが「連休」または「3 連休」の看護師に対してはペナルティの重みを大きくする．

$$f_{\text{n},31}(\boldsymbol{X}) = \sum_{n \in \mathcal{N}_\text{n}} w^{(\text{n})}_{31,n} \sum_{d \in \mathcal{D}_\text{3n}} \max\bigl(0, \sum_{h=0}^{8} x_{n,d+h,\,\text{入}} - 2\bigr)$$

$$w^{(\text{n})}_{31,n} = \begin{cases} 5 & \text{if } f^{(\text{p})}_n = \text{連休 or } f^{(\text{p})}_n = 3\text{ 連休} \\ 1 & \text{otherwise} \end{cases}$$

**S$_\text{n}$-N-32:** 4 連続夜勤の回避

$f_{\text{n},31}$ と同様に，12 時間夜勤に従事する夜勤可能看護師に対し，連続する 13 日間で 4 回の夜勤が割り当てられた場合にペナルティを与える．一方で，16 時間夜勤（「入–明–休」または「入–明–入–明–休」）に従事する夜勤専従看護師に対し，連続する 10 日間で 4 回の夜勤が割り当てられた場合にペナルティを与える．

$$f_{\text{n},32}(\boldsymbol{X}) = \sum_{n \in \mathcal{N}_\text{n}} \sum_{d \in \mathcal{D}_\text{4n}} \max\bigl(0, \sum_{h=0}^{12} x_{n,d+h,\,\text{入}} - 3\bigr) + \sum_{n \in \mathcal{N}_\text{no}} \sum_{d \in \mathcal{D}_\text{4no}} \max\bigl(0, \sum_{h=0}^{9} x_{n,d+h,\,\text{入}} - 3\bigr)$$

**S$_\text{n}$-N-33:** 2 週間で 4 回以上の夜勤の回避

12 時間夜勤（「12h–入–明–休」）に従事する看護師に対し，連続する 14 日間で 4 回の夜勤が割り当てられた場合にペナルティを与える．

$$f_{\text{n},33}(\boldsymbol{X}) = \sum_{n \in \mathcal{N}_\text{n}} \sum_{d \in \mathcal{D}_\text{2w4n}} \max\left(0, \sum_{h=0}^{13} x_{n,d+h,\text{入}} - 3\right)$$

**S$_\text{n}$-N-34:** 夜勤回数に関する違反度の最大最小化

夜勤に関する負担感を公平化するために，12 時間夜勤に従事する看護師の $f_{\text{n},26}$ と $f_{\text{n},28}$ のペナルティ値の和の最大値を別途ペナルティとして与える．

$$f_{\text{n},34}(\boldsymbol{X}) = \max_{n \in \mathcal{N}_\text{n}} \left(f_{\text{n},26}(\boldsymbol{X}_n) + f_{\text{n},28}(\boldsymbol{X}_n)\right)$$

ただし，$\boldsymbol{X}_n$ は，$\boldsymbol{X}$ の要素のうち，$n$ 以外の看護師 $n'\,(\neq n)$ に関する部分を 0 に固定したものである．

**S$_\text{n}$-N-35:** ナース制約に関する違反度の最大最小化

$f_{\text{n},34}$ と同様に，ナース制約に関する違反度を看護師間で公平化するために，$f_{\text{n},26}$ から $f_{\text{n},31}$，および $f_{\text{n},33}$ のペナルティ値の和の最大値を別途ペナルティとして与える．

$$f_{\text{n},35}(\boldsymbol{X}) = \max_{n \in \mathcal{N}_\text{n}} \sum_{i=26, i \neq 32}^{33} f_{\text{n},i}(\boldsymbol{X}_n)$$

## C. 日勤スケジューリングにおけるハード制約の詳細

H$_\text{d}$-N-3, H$_\text{d}$-N-7, H$_\text{d}$-S-8, H$_\text{d}$-S-9 については，夜勤スケジューリングと同様であるため省略する．また，禁止シフト集合 $\mathcal{S}_\text{hf}$ の要素は，夜勤スケジューリングと日勤スケジューリングで異なることに注意する．

**H$_\text{d}$-N-1:** 勤務希望の固定および日勤/休みシフトのみに限定

H$_\text{n}$-N-1 と同様に，勤務希望を固定する．また，日勤スケジューリングでは「日，休」のみを配置するため，それら以外のシフトは配置しないように固定する．

$$\begin{cases} x_{n,d,s} = 1 & \text{if } x'_{n,d,s} = 1 \\ x_{n,d,s} = 0 & \text{if } x'_{n,d,s} = 0, s \notin \{\,\text{日},\text{休}\,\} \end{cases} \quad (n \in \mathcal{N},\ d \in \mathcal{D},\ s \in \mathcal{S})$$

ただし，病院によっては「早，遅」も配置する場合があることに注意する．

**H$_\text{d}$-N-2:** 夜勤専従看護師の日勤制限

夜勤専従看護師に「日」を割り当てない．

$$x_{n,d,\text{日}} = 0 \quad (n \in \mathcal{N}_\text{no},\ d \in \mathcal{D}_\text{t})$$

**H$_\text{d}$-N-4:** 6 連勤の回避

勤務シフトが 6 日連続で割り当てられないようにする．

$$\sum_{h=0}^{5} x_{n,d+h,s} \leq 5 \quad (n \in \mathcal{N},\ d \in \mathcal{D}_\text{6d},\ s \in \mathcal{S}_\text{work})$$

**H$_\text{d}$-N-5:** 次月跨ぎの夜勤を含む 6 連勤の回避

対象月の最終日に「入」が割り当てられており，その前に連続して 4 日以上の日勤が割り当てられた場合，次月の 1 日が「明」であることから次月に跨がる 6 日連続勤務が発生するため，これを回避する．

$$\sum_{h=0}^{3} \sum_{s \in \mathcal{S}_\text{work}} x_{n,d_\text{last}-4+h,s} + x_{n,d_\text{last},\text{入}} \leq 4 \quad (n \in \mathcal{N}_\text{n})$$

**H$_\text{d}$-N-6:** 12 時間日勤の前の 3 連続日勤の回避

12 時間夜勤（「12h–入–明–休」）の場合，12 時間日勤の前に 3 連続日勤が割り当てられた場合，6 連勤となり H$_\text{d}$-N-4 に違反する．しかし，運用上の都合で「12h–休–日–入–明–休」のような例外的なシフトが割り当てられることがある．この場合，12 時間日勤の前に 3 連続日勤が割り当てられたとしても H$_\text{d}$-N-4 によって違反とはみなされないが，非常に負担の大きい勤務となるためこれを回避する．

$$\sum_{h=0}^{2} \sum_{s \in \mathcal{S}_\text{day}} x_{n,d+h,s} + x_{n,d+3,\text{12h}} \leq 3 \quad (n \in \mathcal{N}_\text{n},\ d \in \mathcal{D}_\text{3d})$$

## D. 日勤スケジューリングにおけるソフト制約の詳細

$S_d$-S-1 から $S_d$-S-20, $S_d$-S-22, $S_d$-S-23 は，夜勤スケジューリングにおける定式化と同様であるため省略する．ただし，$S_d$-S-22 におけるソフト禁止シフト集合 $\mathcal{S}_{sf}$ の要素は，夜勤スケジューリングと日勤スケジューリングで異なることに注意する．

**$S_d$-S-21:** リーダー業務可能な看護師の勤務人数の下限

リーダー業務可能な看護師を一人以上配置できない場合にペナルティを与える．土日祝はそもそもの勤務人数が少ないため，リーダー業務可能な看護師を「日，12h」に割り当てることができない場合にペナルティを与える．

$$f_{d,21}(\boldsymbol{X}) = \sum_{d \in \mathcal{D}_t \setminus \mathcal{D}_{t,off}} \max\bigl(0, 1 - \sum_{n \in \mathcal{N}_{dayL}} x_{n,d,日}\bigr) + \sum_{d \in \mathcal{D}_{t,off}} \max\bigl(0, 1 - \sum_{n \in \mathcal{N}_{dayL}} \sum_{s \in \{日,12h\}} x_{n,d,s}\bigr)$$

**$S_d$-N-24:** 週休回数の下限

対象月の週休回数 $N_{off}$ を下回る場合に，その差分に応じてペナルティを与える．

$$f_{d,24}(\boldsymbol{X}) = \sum_{n \in \mathcal{N}} \max\bigl(0, N_{off} - \sum_{d \in \mathcal{D}_t} x_{n,d,休}\bigr)$$

ただし，病院施設によっては，雇用契約上，個人ごとに月の週休回数が異なる場合があるため，それに応じて週休回数を設定する必要があることに注意する．

**$S_d$-N-25:** 週休回数の上限

対象月の週休回数 $N_{off}$ を上回る場合に，その差分に応じてペナルティを与える．

$$f_{d,25}(\boldsymbol{X}) = \sum_{n \in \mathcal{N}} \max\bigl(0, \sum_{d \in \mathcal{D}_t} x_{n,d,休} - N_{off}\bigr)$$

**$S_d$-N-26:** 土日祝における休み回数の下限

家族等との時間を確保をするため，土日祝に割り当てられた「休，特休」の合計回数が 4 回未満となる場合に，その差分に応じてペナルティを与える．

$$f_{d,26}(\boldsymbol{X}) = \sum_{n \in \mathcal{N}} \max\bigl(0, 4 - \sum_{d \in \mathcal{D}_{t,off}} \sum_{s \in \mathcal{S}_{off}} x_{n,d,s}\bigr)$$

ここで，月当たりに最低 8 日の土日祝が存在するため 4 回としている．

**$S_d$-N-27:** 3 連続日勤の回避

日勤帯は手術や検査等のイベントが多く，勤務強度が高いため，3 日以上日勤帯に当たる勤務「日，12h，早，遅」が連続すると非常に負担が大きい．そこで，それらが 3 連続で割り当てられた場合にペナルティを与える．

$$f_{d,27}(\boldsymbol{X}) = \sum_{n \in \mathcal{N}} \sum_{d \in \mathcal{D}_{3d}^{(d)}} \max\bigl(0, \sum_{h=0}^{2} \sum_{s \in \mathcal{S}_{day} \setminus \{他\}} x_{n,d+h,s} - 2\bigr)$$

本システムでは，時短看護師（月から金での固定勤務）や休みの好みが「7 日 2 休」である日勤専従看護師を対象外とし，余分にペナルティを与えないようにしている．

**$S_d$-N-28:** 5 連勤の後の単発休みの回避

5 連勤の後に単発休みしか割り当てられない場合，十分な休息が取れない．そのため，このような勤務が割り当てられた場合にペナルティを与える．

$$f_{d,28}(\boldsymbol{X}) = \sum_{n \in \mathcal{N}} \sum_{d \in \mathcal{D}_{5d1off}} \max\bigl(0, \sum_{h=0}^{4} \sum_{s \in \mathcal{S}_{work}} x_{n,d+h,s} + \sum_{s \in \mathcal{S}_{off}} x_{n,d+5,s} + \sum_{s \in \mathcal{S}_{work}} x_{n,d+6,s} - 6\bigr)$$

**$S_d$-N-29:** 9 日間における週休回数の下限

「日–12h–入–明–休–日–12h–入–明–休」のようなシフトでは，9 日間において休み回数が 1 回となるため，十分な休息が取れない．そのため，9 日間における休み回数が 2 回未満の場合にペナルティを与える．

$$f_{d,29}(\boldsymbol{X}) = \sum_{n \in \mathcal{N}} \sum_{d \in \mathcal{D}_{9d2off}} \max\bigl(0, 2 - \sum_{h=0}^{8} \sum_{s \in \mathcal{S}_{off}} x_{n,d+h,s}\bigr)$$

**S$_\text{d}$-N-30:** 土日連休回数の下限

土日に連続して休みが割り当てられているか否かを表す二値変数 $s_{n,d}^{(\text{ss})}$ の合計が 1 未満の場合にペナルティを与える．このとき，休みの好みが「土日連休」である看護師に対しては，ペナルティの重みを大きくする．

$$f_{\text{d},30}(\boldsymbol{X}) = \sum_{n \in \mathcal{N}} w_{30,n}^{(\text{d})} \max\bigl(0, 1 - \sum_{d \in \mathcal{D}_\text{sat}} s_{n,d}^{(\text{ss})}\bigr)$$

$$w_{30,n}^{(\text{d})} = \begin{cases} 5 & \text{if } f_n^{(\text{p})} = \text{土日連休} \\ 1 & \text{otherwise} \end{cases}$$

ただし，$n \in \mathcal{N}$, $d \in \mathcal{D}_\text{sat}$ に対し，以下を満たす．

$$\sum_{s \in \mathcal{S}_\text{off}} x_{n,d,s} + \sum_{s \in \mathcal{S}_\text{off}} x_{n,d+1,s} \geq 2 \cdot s_{n,d}^{(\text{ss})}, \quad \sum_{s \in \mathcal{S}_\text{off}} x_{n,d,s} + \sum_{s \in \mathcal{S}_\text{off}} x_{n,d+1,s} \leq 1 + s_{n,d}^{(\text{ss})}$$

**S$_\text{d}$-N-31:** 2 連休回数の下限

2 連休が割り当てられているか否かを表す二値変数 $s_{n,d}^{(\text{2off})}$ の合計が 2 未満の場合，その差分をペナルティとして与える．このとき，休みの好みが「2 連休」である看護師に対しては，ペナルティの重みを大きくする．

$$f_{\text{d},31}(\boldsymbol{X}) = \sum_{n \in \mathcal{N}} w_{31,n}^{(\text{d})} \max\bigl(0, 2 - \sum_{d \in \mathcal{D}_\text{2off}^{(\text{d})}} s_{n,d}^{(\text{2off})}\bigr)$$

$$w_{31,n}^{(\text{d})} = \begin{cases} 5 & \text{if } f_n^{(\text{p})} = 2\text{ 連休} \\ 1 & \text{otherwise} \end{cases}$$

ただし，$n \in \mathcal{N}$, $d \in \mathcal{D}_\text{2off}^{(\text{d})}$ に対し，以下を満たす．

$$\begin{cases} \displaystyle\sum_{s \in \mathcal{S}_\text{work}} x_{n,d,s} + \sum_{h=1}^{2} \sum_{s \in \mathcal{S}_\text{off}} x_{n,d+h,s} + \sum_{s \in \mathcal{S}_\text{work}} x_{n,d+3,s} \geq 4 \cdot s_{n,d}^{(\text{2off})} \\ \displaystyle\sum_{s \in \mathcal{S}_\text{work}} x_{n,d,s} + \sum_{h=1}^{2} \sum_{s \in \mathcal{S}_\text{off}} x_{n,d+h,s} + \sum_{s \in \mathcal{S}_\text{work}} x_{n,d+3,s} \leq 3 + s_{n,d}^{(\text{2off})} \end{cases}$$

**S$_\text{d}$-N-32:** 3 連続以上の休み回数の下限

3 連続以上の休みが割り当てられているか否かを表す二値変数 $s_{n,d}^{(\text{3off})}$ の合計が 1 未満の場合，その差分をペナルティとして与える．このとき，休みの好みが「3 連休」である看護師に対しては，ペナルティの重みを大きくする．

$$f_{\text{d},32}(\boldsymbol{X}) = \sum_{n \in \mathcal{N}} w_{32,n}^{(\text{d})} \max\bigl(0, 1 - \sum_{d \in \mathcal{D}_\text{3off}^{(\text{d})}} s_{n,d}^{(\text{3off})}\bigr)$$

$$w_{32,n}^{(\text{d})} = \begin{cases} 5 & \text{if } f_n^{(\text{p})} = 3\text{ 連休} \\ 1 & \text{otherwise} \end{cases}$$

ただし，$n \in \mathcal{N}$, $d \in \mathcal{D}_\text{3off}^{(\text{d})}$ に対し，以下を満たす．

$$\sum_{h=0}^{2} \sum_{s \in \mathcal{S}_\text{off}} x_{n,d+h,s} \geq 3 \cdot s_{n,d}^{(\text{3off})}, \quad \sum_{h=0}^{2} \sum_{s \in \mathcal{S}_\text{off}} x_{n,d+h,s} \leq 2 + s_{n,d}^{(\text{3off})}$$

**S$_\text{d}$-N-33:** 単発日勤回数の上限

「休–日–休」のような単発日勤の数が増加すると，その分だけ連休回数が減少するため，単発日勤か否かを表す二値変数 $s_{n,d}^{(\text{1d})}$ の合計が，休みの好みに基づいて設定された上限回数を上回る場合，その差分をペナルティとして与える．休みの好みが「単発休み」以外の場合は上限を 0 とし，「単発休み」の場合は上限を 2 とした上でペナルティの重みも大きくする．

$$f_{\text{d},33}(\boldsymbol{X}) = \sum_{n \in \mathcal{N}_\text{so}} w_{33,n}^{(\text{d})} \max\bigl(0, \sum_{d \in \mathcal{D}_\text{1d}} s_{n,d}^{(\text{1day})} - 2\bigr) + \sum_{n \in \mathcal{N} \setminus \mathcal{N}_\text{so}} w_{33,n}^{(\text{d})} \sum_{d \in \mathcal{D}_\text{1d}} s_{n,d}^{(\text{1day})}$$

$$w_{33,n}^{(\text{d})} = \begin{cases} 5 & \text{if } f_n^{(\text{p})} = \text{単発休み} \\ 1 & \text{otherwise} \end{cases}$$

ただし，$n \in \mathcal{N}$, $d \in \mathcal{D}_{\text{1d}}$ に対し，以下を満たす．

$$\begin{cases} \displaystyle\sum_{s \in \mathcal{S}_{\text{off}}} x_{n,d,s} + \sum_{s \in \mathcal{S}_{\text{day}}} x_{n,d+1,s} + \sum_{s \in \mathcal{S}_{\text{off}}} x_{n,d+2,s} \geq 3 \cdot s_{n,d}^{(\text{1day})} \\ \displaystyle\sum_{s \in \mathcal{S}_{\text{off}}} x_{n,d,s} + \sum_{s \in \mathcal{S}_{\text{day}}} x_{n,d+1,s} + \sum_{s \in \mathcal{S}_{\text{off}}} x_{n,d+2,s} \leq 2 + s_{n,d}^{(\text{1day})} \end{cases}$$

**S$_{\text{d}}$-N-34:** ナース制約に関する違反度の最大最小化

S$_{\text{n}}$-N-35 と同様に，ナース制約に関する違反度を看護師間で公平化するために，$f_{\text{d},26}$ から $f_{\text{d},28}$ および $f_{\text{d},30}$ から $f_{\text{d},33}$ のペナルティ値の和の最大値を別途ペナルティとして与える．

$$f_{\text{d},34}(\boldsymbol{X}) = \max_{n \in \mathcal{N}} \sum_{i=26, i \neq 29}^{33} f_{\text{d},i}(\boldsymbol{X}_n)$$